\input amstex
\NoBlackBoxes
\pageno=1

\loadbold
\documentstyle{amsppt}
\def\varep{\varepsilon}
\def\eps{\varepsilon}

\def\e{{\text{\rm e}}}
\def\trd{{\text{\rm d}}}
\def\trc{{\text{\rm c}}}
\def\loc{_{\text{\rm loc}}}
\def\sign{\operatorname{sgn}}
\def\capt{\operatorname{cap}}
\def\meas{\operatorname{meas}}
\def\supp{\operatorname{supp}}

\def\ds{\displaystyle}
\def\bs{\boldsymbol}
\def\nd{\noindent}

\topmatter
\title{Nonlinear elliptic equations with measures revisited
}\endtitle
\author{Ha\"im Brezis$^{(1),(2)}$, Moshe Marcus$^{(3)}$ and Augusto
C. Ponce$^{(4)}$}\endauthor
\medskip
\medskip

\address
${}^{\text {(1)}}$  Rutgers University\endgraf
Dept. of Mathematics, Hill Center, Busch Campus\endgraf
110 Frelinghuysen Rd., Piscataway, NJ 08854, USA\endgraf
\endaddress
\email
brezis\@math.rutgers.edu\endemail
\null

\address
${}^{\text {(2)}}$  Laboratoire J.-L. Lions\endgraf
Universit\'e P. et M. Curie, B.C. 187\endgraf
4 Pl. Jussieu\endgraf
75252 Paris Cedex 05, France\endgraf
\endaddress
\email
brezis\@ccr.jussieu.fr
\endemail

\address ${}^{\text{(3)}}$ Technion\endgraf
Dept. of Mathematics\endgraf
Haifa 32000, Israel\endgraf
\endaddress
\email
marcusm\@tx.technion.ac.il\endemail

\address ${}^{\text{(4)}}$
Institute for Advanced Study\endgraf
Princeton, NJ 08540, USA\endgraf
\endaddress
\email
augponce\@math.ias.edu\endemail
\null

\abstract
We study the existence of solutions of the nonlinear problem
$$
\left\{
\alignedat2
-\Delta u + g(u)  & = \mu & &\quad \text{in }  \Omega,\\
\qquad\qquad\;\;\,  u & = 0 & & \quad \text{on } \partial \Omega,
\endaligned
\right.\tag P
$$
where $\mu$ is a Radon measure and $g : \Bbb R \to \Bbb R$ is a
nondecreasing continuous function with $g(0) = 0$. This equation need
not have a solution for every measure $\mu$, and we say that $\mu$ is
a good measure if (P) admits a solution. We show that for every $\mu$
there exists a largest good measure $\mu^* \leq \mu$. This {\it
reduced measure}\/ has a number of remarkable properties.

To cite this paper: Ha\"\i m Brezis, Moshe Marcus and Augusto C.~Ponce,
{\it Nonlinear elliptic equations with measures
  revisited}. In: Mathematical Aspects of Nonlinear Dispersive Equations (J.~Bourgain, C.~Kenig, S.~Klainerman, eds.), Annals of Mathematics Studies, 163, Princeton University Press, Princeton, NJ, 2007, pp.~55--110.
\endabstract

\toc
\head{} 0. Introduction \endhead
\head{} 1. Construction of $u^*$ and $\mu^*$.  Proofs of Proposition 1 and Theorems 1, 2 \endhead
\head{} 2. Good measures. Proofs of Theorems 4, 6 \endhead
\head{} 3. Some properties of the mapping $\mu\mapsto \mu^*$ \endhead
\head{} 4. Approximation via $\rho_n * \mu$ \endhead
\head{} 5. Further convergence results \endhead
\head{} 6. Measures which are good for every $g$ must be diffuse \endhead
\head{} 7. Signed measures and general nonlinearities $g$ \endhead
\head{} 8. Examples \endhead
\head{} 9. Further directions and open problems \endhead
\head{} Appendix A: Decomposition of measures into diffuse and
concentrated parts \endhead
\head{} Appendix B: Standard existence, uniqueness and comparison
results \endhead
\head{} Appendix C: Correspondence between $\big[ C_0(\overline\Omega)\big]^*$
and $\big[ C(\overline\Omega) \big]^*$ \endhead
\head{} Appendix D: A new decomposition for diffuse measures \endhead
\head{} Appendix E: Equivalence between $\capt_{H^1}$ and
$\capt_{\Delta,1}$ \endhead
\head{} References \endhead
\endtoc

\endtopmatter
\document


\subhead 0. Introduction\endsubhead
\medskip

Let $\Omega\subset \Bbb R^N$ be a bounded domain with smooth boundary.
Let $g: \Bbb R \to
\Bbb R$ be a continuous, nondecreasing function such that $g(0)= 0$.
In this paper we are concerned with the problem
$$
\left\{
\alignedat2
-\Delta u + g(u)  & = \mu & &\quad \text{in }  \Omega,\\
\qquad\qquad\;\;\,  u & = 0 & & \quad \text{on } \partial \Omega,
\endaligned
\right.\tag 0.1
$$
where $\mu$ is a measure.  The study of (0.1) when
$\mu \in L^1(\Omega)$ was initiated by Brezis-Strauss [BS]; their main
result asserts that for {\it every} $\mu \in L^1$ and {\it
every} $g$ as above, problem (0.1) admits a unique weak solution (see
Theorem~B.2 in Appendix~B below).  The
right concept of weak solution is the following:
$$
\cases
u \in L^1(\Omega),\; g(u) \in L^1(\Omega) \text{ and } \\
\ds  -\int_\Omega u \Delta \zeta + \int_\Omega g(u) \zeta = \int_\Omega
\zeta \, d\mu \quad \forall \zeta \in C^2(\overline\Omega),\; \zeta = 0 \text
{ on } \partial \Omega.
\endcases\tag  0.2
$$
It will be convenient to write
$$
C_0(\overline\Omega) = \big\{ \zeta \in C (\overline\Omega) \; ; \, \zeta = 0
\text{ on } \partial \Omega \big\}
$$
and
$$
C_0^2 (\overline\Omega) = \big\{ \zeta \in C^2 (\overline\Omega) \; ; \, \zeta = 0
\text{ on } \partial \Omega \big\},
$$
and to say that (0.1) holds in the sense of $(C^2_0)^*$.  We
will often omit the word ``weak'' and simply say that $u$ is a solution of
(0.1), meaning (0.2).  It follows from standard (linear) regularity theory that a weak solution $u$ belongs to $W^{1,
q}_0(\Omega)$ for every $q< \frac{N}{N-1}$ (see, e.g., [S] and Theorem~B.1 below).

The case where $\mu$ is a measure turns out to be much more subtle
than one might expect.  It was observed in 1975 by
Ph.~B\'enilan and H.~Brezis (see [B1], [B2], [B3], [B4], [BB]
and Theorem~B.6 below) that if $N\geq 3$ and $g(t) = |t|^{p-1} t$ with
$p\geq \frac{N}{N-2}$, then (0.1) has {\it no solution}\/ when $\mu = \delta_a$, a
Dirac mass at a point $a\in \Omega$.  On the other hand, it was also
proved (see Theorem~B.5 below) that if $g(t) = |t|^{p-1} t$ with $p< \frac{N}{N-2}$
(and $N\geq 2$), then (0.1) has a solution for any measure $\mu$.
Later Baras-Pierre [BP] (see also [GM]) {\it characterized}\/ all
measures $\mu$ for which
(0.1) admits a solution.  Their necessary and sufficient condition for
the existence of a solution
when $p\geq \frac{N}{N-2}$ can be expressed in two equivalent ways:
$$
\cases
\mu \text{ admits a decomposition } \mu = f_0 - \Delta v_0
\text{ in the } (C^2_0)^* \text{-sense},\\ \text{with } f_0\in L^1 \text{ and }  v_0\in L^p,
\endcases\tag 0.3
$$
or
$$
|\mu| (A) = 0\quad \text{for every Borel set $A \subset \Omega$ with $\capt_{2,p'}{(A)} =
 0$,}\tag 0.4
$$
where $\capt_{2,p'}$ denotes the capacity associated to $W^{2,p'}$.

Our goal in this paper is to analyze the nonexistence
mechanism and to describe what happens if one ``forces'' (0.1) to
have a solution in cases where the equation ``refuses'' to possess
one.  The natural approach is to introduce an approximation scheme.
For example, $\mu$ is kept fixed and $g$ is truncated.  Alternatively,
$g$ is kept fixed and $\mu$ is approximated, e.g., via convolution.  It
was originally observed by one of us (see [B4]) that if $N\geq 3$, $g(t) =
|t|^{p-1} t$, with $p\geq \frac{N}{N-2}$, and $\mu = \delta_a$, with $a
\in \Omega$, then all ``natural" approximations $(u_n)$ of (0.1) converge
to $u \equiv 0$.  And, of course, $u\equiv 0$ is {\it not}\/ a solution
of (0.1) corresponding to $\mu = \delta_a$~! It is this kind of
phenomenon that we propose to explore in full generality. We are led
to study the convergence of the approximate solutions $(u_n)$ under
various assumptions on the sequence of data.

Concerning the function $g$ we will assume {\it throughout the rest of
the paper}\/ (except in Section~7) that $g: \Bbb R\to \Bbb R$ is
continuous, nondecreasing, and
that
$$
g(t) = 0\quad \forall t \leq 0.\tag 0.5
$$

\remark{Remark 1}  Assumption (0.5) is harmless when the data $\mu
$ is nonnegative, since the corresponding solution $u$ is nonnegative by the maximum
principle and it is only the restriction of $g$ to $[0,\infty)$ which
is relevant.  However when $\mu$ is a signed measure it is worthwhile
to remove assumption (0.5) and this is done in Section~7 below.
\endremark

\medskip

By a {\it measure}\/ $\mu$ we mean a continuous linear functional on
$C_0(\overline\Omega),$ or equivalently a finite measure on
$\overline\Omega$ such that $|\mu|(\partial \Omega) = 0$ (see
Appendix~C below).  The space
of measures is denoted by $\Cal M(\Omega)$ and is equipped with the standard norm
$$
\|\mu\|_{\Cal M} = \sup{\bigg \{\int_\Omega \varphi \, d \mu\, ;\, \varphi \in
C_0(\overline\Omega) \text{ and } \|\varphi\|_{L^\infty} \leq 1\bigg \}}.
$$

By a (weak) {\it solution}\/ $u$ of (0.1) we mean that (0.2) holds.  A
(weak) {\it subsolution}\/ $u$ of (0.1) is a function $u$ satisfying
$$
\cases
u \in L^1(\Omega), \; g(u) \in L^1(\Omega) \text{ and } \\
\ds -\int_\Omega u \Delta \zeta + \int_\Omega g(u) \zeta \leq \int_\Omega
\zeta \, d\mu \quad\forall \zeta \in C^2_0(\overline\Omega),
\; \zeta \geq 0 \text{ in } \Omega.
\endcases\tag 0.6
$$

We will say that $\mu \in \Cal M(\Omega)$ is a {\it good measure} if
(0.1) admits a solution.  If $\mu$ is a good measure, then equation
(0.1) has exactly one solution $u$ (see
Corollary~B.1 in Appendix B). We denote by $\Cal G$ the set of good
measures (relative to $g$).

\remark{Remark 2}
In many places throughout this paper, the quantity $\int_\Omega \zeta
\, d\mu$, with $\zeta \in C_0^2(\overline\Omega)$, plays an important
role. Such an expression makes sense even for measures $\mu$ which are
not bounded but merely locally bounded in $\Omega$, and such that
$\int_\Omega \rho_0 \, d|\mu| < \infty$, where $\rho_0(x) = d(x ,
\partial\Omega)$. Many of our
results remain valid for such measures provided some of the
statements (and the proofs) are slightly modified. In this case, the
condition $g(u) \in L^1(\Omega)$ in (0.2) (and also in (0.6)) must be
replaced by $g(u) \rho_0 \in L^1(\Omega)$. Since we have not pursued
this direction, we shall leave the details to the reader.
\endremark

\medskip
In Section 1 we will introduce the first approximation method, namely
$\mu$ is fixed and $g$ is ``truncated''.  In the sequel we denote by
$(g_n)$ a sequence of functions $g_n: \Bbb R\to \Bbb R$ which are
continuous, nondecreasing and satisfy the following conditions:
$$
\gather
0\leq g_1(t) \leq g_2(t)\leq \ldots\leq g(t)\quad\forall t \in \Bbb R,
\tag 0.7\\
g_n(t) \to g(t)\quad \forall t \in \Bbb R. \tag 0.8
\endgather
$$
(Recall that, by Dini's lemma, conditions (0.7) and (0.8) imply that
$g_n \to g$ uniformly on compact subsets of $\Bbb R$).

\noindent
If $N \geq 2$, we assume in addition that each $g_n$ has subcritical
growth, i.e., that there exist $C>0$ and  $p<\frac{N}{N-2}$ (possibly
depending on $n$) such that
$$
g_n(t)\leq C(|t|^p+ 1)\quad \forall t \in \Bbb R. \tag 0.9
$$

\medskip
\noindent
A good example to keep in mind is $g_n(t) = \min{\{ g(t), n\}}$,
$\forall t \in \Bbb R$.

\bigskip

Our first result is

\proclaim{Proposition 1}  Given any measure $\mu \in \Cal M(\Omega)$,
let $u_n$ be the unique solution of
$$
\left\{
\alignedat2
-\Delta u_n + g_n(u_n) & = \mu & & \quad \text{in } \Omega,\\
 u_n & = 0 & & \quad               \text{on } \partial\Omega.
\endalignedat
\right.\tag 0.10
$$
Then $u_n\downarrow u^*$ in $\Omega$ as $n\uparrow \infty$, where
$u^*$ is the largest subsolution of\/ $(0.1)$.  Moreover we have
$$
\Big|\int_\Omega u^*\Delta \zeta\Big| \leq 2 \|\mu\|_{\Cal
M}\|\zeta\|_{L^\infty}\quad \forall \zeta \in C^2_0
(\overline\Omega)\tag0.11
$$
and
$$
\int_\Omega g(u^*) \leq \|\mu\|_{\Cal M}.\tag 0.12
$$
\endproclaim

An important consequence of Proposition 1 is that $u^*$ {\it does not
depend on the choice of the truncating sequence $(g_n)$.  It is an
intrinsic object which will play an important role in the sequel}. In
some sense, $u^*$ is the ``best one can do"~(!) in the absence of a
solution.

\remark{Remark 3}  If $\mu$ is a good measure, then $u^*$ coincides
with the unique solution $u$ of (0.1); this is an easy consequence of
standard comparison arguments (see Corollary~B.2 in Appendix~B).
\endremark

\medskip
We now introduce the basic concept of {\it reduced measure}. From
(0.11), (0.12), and the density of $C_0^2(\overline\Omega)$ in
$C_0(\overline \Omega)$ (easy to check), we see that there exists a unique
measure $\mu^* \in \Cal M (\Omega)$ such that
$$
-\int_\Omega u^* \Delta \zeta + \int_\Omega g(u^*)\zeta = \int_\Omega
 \zeta \, d \mu^*\quad \forall \zeta \in C^2_0 (\overline\Omega).\tag0.13
$$

We call $\mu^*$ the reduced measure associated to $\mu$.
Clearly, $\mu^*$ is always a good measure.  Since $u^*$ is a
subsolution of (0.1), we have
$$
\mu^*\leq \mu.\tag 0.14
$$
Even though we have not indicated the dependence on $g$ we emphasize
that $\mu^*$ {\it does depend} on $g$ (see Section~8 below).
\medskip

One of our main results is

\proclaim{Theorem 1}  The reduced measure $\mu^*$ is the largest good measure
$\leq \mu$.
\endproclaim

Here is an easy consequence:

\proclaim{Corollary 1}  We have
$$
0\leq \mu-\mu^* \leq \mu^+ = \sup{\{ \mu, 0\}}.\tag 0.15
$$
In particular,
$$
|\mu^*|\leq |\mu|\tag 0.16
$$
and
$$
[\mu \geq 0 ] \quad \Longrightarrow \quad [\mu^*
\geq 0]. \tag 0.17
$$
\endproclaim

\smallskip
Indeed, every measure $\nu \leq 0$ is a good measure since the
solution $v$ of
$$
\left\{
\alignedat2
-\Delta v & = \nu & & \quad \text{in } \Omega,\\
v & = 0 & & \quad \text{on } \partial \Omega,
\endalignedat
\right.
$$
satisfies $v\leq 0 \text{ in } \Omega$, and therefore by (0.5)
$$
-\Delta v + g(v) = \nu \quad \text{in } (C_0^2)^*.
$$
In particular, $-\mu^-$ is a good measure (recall that $ \mu^- =
\sup{\{ -\mu, 0\}}$). Since $-\mu^- \leq \mu$, we deduce from Theorem 1 that
$$
-\mu^- \leq \mu^*,
$$
and consequently
$$
\mu - \mu^* \leq \mu + \mu^- = \mu^+.
$$

\medskip
Our next result asserts that the measure $\mu - \mu^*$ is concentrated on a
small set:

\proclaim{Theorem 2}  There exists a Borel set $\Sigma \subset \Omega$
with $\capt{(\Sigma)} = 0$ such that
$$
(\mu- \mu^*)(\Omega\setminus \Sigma) = 0.\tag 0.18
$$
\endproclaim
Here and throughout the rest of the paper ``cap" denotes the Newtonian
($H^1$) capacity with respect to $\Omega$.
\medskip

\remark{Remark 4}
Theorem~2 is optimal in the following sense. Given any measure $\mu
\geq 0$ concentrated on a set of zero capacity, there exists some $g$
such that $\mu^* = 0$ (see Theorem~14 below). In particular, $\mu -
\mu^*$ can be {\it any}\/ nonnegative measure concentrated on a set of zero capacity.
\endremark

\medskip

Here is a useful

\proclaim{Definition}  A measure $\mu\in \Cal M(\Omega)$ is called
diffuse if\/ $|\mu|(A) = 0$ for every Borel set $A\subset\Omega$  such
that\/ $\capt{(A)} =0$.
\endproclaim

An immediate consequence of Corollary 1 and Theorem 2 is

\proclaim{Corollary 2}  Every diffuse measure $\mu \in \Cal M(\Omega)$
is a good measure.
\endproclaim

Indeed, let $\Sigma$ be as in Theorem 2, so that $\capt{(\Sigma)}=0$ and
$$
(\mu - \mu^*)(\Omega\setminus \Sigma)=0.
$$
On the other hand, (0.15) implies
$$
(\mu - \mu^*)(\Sigma)\leq \mu^+ (\Sigma) = 0,
$$
since $\mu$ is diffuse.  Therefore
$$
(\mu - \mu^*)(\Omega) = 0,
$$
so that $\mu = \mu^*$ and thus $\mu$ is a good measure.

\remark{Remark 5}  The converse of Corollary 2 is {\it not}\/ true.  In
Example 5 (see Section~8 below) the measure $\mu = c\delta_a$, with $0< c \leq
4\pi$ and $a\in \Omega$, is a good measure, but it is not diffuse
--- $\capt{(\{a\})}= 0$, while $\mu(\{a\}) = c > 0$. See, however,
Theorem~5.
\endremark

\remark{Remark 6}  Recall that a measure $\mu$ is diffuse if and only
if $\mu \in L^1 + H^{-1}$; more precisely, there exist $f_0\in
L^1(\Omega)$ and $v_0\in H_0^1(\Omega)$ such that
$$
\int_\Omega \zeta \, d \mu = \int_\Omega f_0 \zeta - \int_\Omega
\nabla v_0\cdot
\nabla \zeta\quad\forall \zeta \in C_0(\overline \Omega) \cap H^1_0.\tag 0.19
$$
The implication $[\mu \in L^1 + H^{-1}] \Rightarrow [\mu \text{
diffuse}]$ is due to Grun-Rehomme~[GRe].  (In fact he proved only that
$[\nu \in H^{-1}] \Rightarrow [\nu \text{ diffuse}],$ but
$L^1$-functions are diffuse measures --- since $[\capt{(A)}=0]
\Rightarrow [|A| = 0]$ --- and the sum of two diffuse measures is
diffuse).  The converse [$\mu$
diffuse] $\Rightarrow [\mu \in L^1 + H^{-1}]$ is due to
Boccardo-Gallou\"et-Orsina [BGO1] (and was suggested by earlier results
of Baras-Pierre [BP] and Gallou\"et-Morel [GM]).  As a consequence of
Corollary 2 we obtain that, for every measure $\mu$ of the form (0.19),
the problem
$$
\left\{
\alignedat2
-\Delta u + g(u) & = \mu& & \quad  \text{in } \Omega,\\
u & = 0 & & \quad \text{on } \partial\Omega,
\endalignedat
\right. \tag 0.20
$$
admits a unique solution.  In fact, the same conclusion was already
known for any {\it distribution}\/ in $L^1+ H^{-1}$, not
necessarily in $\Cal M(\Omega)$. (The proof, which combines techniques
from Brezis-Browder [BBr] and Brezis-Strauss [BS], is sketched in
Appendix B below; see Theorem~B.4). A very useful sharper version of the
[BGO1] decomposition is the following:

\proclaim{Theorem 3}
Assume $\mu \in \Cal M(\Omega)$ is a diffuse measure. Then, there exist $f
\in L^1(\Omega)$ and $v \in C_0(\overline\Omega) \cap H^1_0$ such that
$$
\int_\Omega \zeta \, d \mu = \int_\Omega f \zeta - \int_\Omega
\nabla v \cdot
\nabla \zeta\quad\forall \zeta \in C_0(\overline \Omega) \cap H^1_0. \tag{0.21}
$$
In addition, given any $\delta >0$, then $f$ and $v$ can be chosen so that
$$
\| f \|_{L^1} \leq \| \mu \|_{\Cal M},  \quad \|v\|_{L^\infty} \leq
\delta \|\mu\|_{\Cal M} \quad \text{and} \quad \|v\|_{H^1} \leq
\delta^{1/2} \|\mu\|_{\Cal M}. \tag 0.22
$$
\endproclaim

The proof of Theorem~3 is presented in Appendix~D below.
\endremark

\smallskip

In Section 2 we present some basic properties of the good measures.
Here is a first one:

\proclaim{Theorem 4}  Suppose $\mu_1$ is a good measure.
Then any measure $\mu_2\leq \mu_1$ is also a good measure.
\endproclaim

We now deduce a number of consequences:

\proclaim{Corollary 3} Let $\mu \in \Cal M(\Omega)$. If $\mu^+$ is
diffuse, then $\mu$ is a good measure.
\endproclaim

In fact, by Corollary~2, $\mu^+$ diffuse implies that $\mu^+$ is a
good measure. Since $\mu \leq \mu^+$, it follows from Theorem~4 that
$\mu$ is a good measure.

\proclaim{Corollary 4}  If $\mu_1$ and $\mu_2$ are good measures, then
so is $\nu = \sup{\{\mu_1, \mu_2\}}$.
\endproclaim

Indeed, by Theorem 1 we have $\mu_1 \leq \nu^*$ and $\mu_2 \leq
\nu^*$.  Thus $\nu\leq \nu^*\leq \nu$,  and hence $\nu = \nu^*$ is
good measure.

\proclaim{Corollary 5} The set $\Cal G$ of good measures is convex.
\endproclaim

Indeed, let $\mu_1, \mu_2\in \Cal G$.  For any $t\in [0,1]$, we have
$$
t\mu_1 + (1 - t) \mu_2 \leq \sup{\{\mu_1, \mu_2\}}.
$$
Applying Corollary~4 and Theorem~4, we deduce that $t\mu_1 + (1-t)\mu_2
\in \Cal G$.

\proclaim{Corollary 6} For every measure $\mu\in \Cal M(\Omega)$ we
have
$$
\|\mu - \mu^*\|_\Cal M = \min_{\nu \in \Cal G}{\|\mu - \nu\|_{\Cal
M}}. \tag 0.23
$$
Moreover, $\mu^*$ is the unique good measure which achieves the minimum.
\endproclaim

\demo{Proof}  Let $\nu \in \Cal G$ and write
$$
|\mu-\nu| = (\mu-\nu)^+ + (\mu-\nu)^- \geq (\mu - \nu)^+ = \mu -
 \inf{\{\mu, \nu\}}.
$$
But $\tilde \nu=\inf{\{\mu, \nu\}}\in \Cal G$ by Theorem~4.  Applying
Theorem 1 we find $\tilde \nu \leq \mu^*$.  Hence
$$
|\mu - \nu|\geq \mu - \tilde \nu \geq \mu- \mu^*\geq 0,
$$
and therefore
$$
\|\mu - \nu\|_{\Cal M} \geq \| \mu - \mu^*\|_\Cal M\, ,
$$
which gives (0.23). In order to establish uniqueness, assume $\nu \in
\Cal G$ attains
the minimum in (0.23). Note that $\inf{ \{\mu,\nu\} } $ is a good
measure $\leq \mu$ and
$$
\big\| \mu - \inf{ \{\mu,\nu\} } \big\|_{\Cal M} \leq \| \mu - \nu
\|_{\Cal M}.
$$
Thus, $\nu = \inf{ \{\mu,\nu\} } \leq \mu$. By Theorem~1, we deduce
that $\nu \leq \mu^* \leq \mu$. Since $\nu$ achieves the minimum in
(0.23), we must have $\nu = \mu^*$.

\enddemo
\medskip

As we have already pointed out, the set $\Cal G$ of good measures
associated to (0.1) depends on the nonlinearity $g$. Sometimes, in order to
emphasize this dependence, we shall denote $\Cal G$ by $\Cal G(g)$. By
Corollary~3, if $\mu \in \Cal M(\Omega)$ and $\mu^+$ is diffuse, then
$\mu \in \Cal G(g)$ for every $g$ satisfying (0.5).
The converse is also true. More precisely,

\proclaim{Theorem 5}
Let $\mu \in \Cal M (\Omega)$. Then $ \mu \in \Cal G(g)$ for every $g$
if and only if $\mu^+$ is diffuse.
\endproclaim

\smallskip

We also have a characterization of good measures in the spirit
of the Baras-Pierre result (0.3):

\proclaim{Theorem 6}  A measure $\mu \in \Cal M(\Omega)$ is a good
measure if and only if $\mu$ admits a decomposition
$$
\mu = f_0 -\Delta v_0 \quad \text{in } \Cal D'(\Omega),
$$
with $f_0\in L^1(\Omega)$, $v_0\in L^1(\Omega)$ and $g(v_0) \in
L^1(\Omega)$.
\endproclaim

\medskip

\proclaim{Corollary 7} We have
$$
\Cal G + L^1(\Omega) \subset \Cal G.
$$
\endproclaim
\medskip

In Section~3 we discuss some properties of the mapping $\mu\mapsto \mu^*$.
For example, we show that for every $\mu, \nu \in \Cal
M(\Omega)$, we have
$$
(\mu^*-\nu^*)^+ \leq (\mu - \nu)^+.
\tag 0.24
$$
Inequality (0.24) implies, in particular, that
$$
[\mu \leq \nu] \quad \Longrightarrow \quad [\mu^*\leq \nu^*]\tag 0.25
$$
and
$$
|\mu^*- \nu^*|\leq |\mu- \nu|.\tag 0.26
$$
\medskip
In Section~4 we examine another approximation scheme. We now
keep $g$ fixed but we smooth $\mu$ via convolution.  Let $\mu_n =
\rho_n *\mu$ and let $u_n$ be the solution of
$$
\left\{
\alignedat2
-\Delta u_n + g(u_n) & = \mu_n & & \quad \text{in } \Omega,\\
  u_n & =0 & & \quad \text{on } \partial \Omega.
\endalignedat
\right. \tag 0.27
$$
We prove (assuming in addition $g$ is convex) that $u_n
\to u^*$ in $L^1(\Omega)$, where $u^*$ is given by Proposition 1. In
Section 5 we discuss other convergence results.

Theorem~5 is established in Section~6.
In Section~7 we extend Proposition~1 to deal with the case where $\mu
\in \Cal M(\Omega)$ is a signed measure, but assumption (0.5) is no
longer satisfied. Finally, in Section 8 we present several examples
where the measure $\mu^*$ can be explicitly identified and in
Section~9 we propose various directions of research.

Part of the results in this paper were announced in [BMP].

\bigskip


\subhead 1. Construction of $\bs u^{\pmb *}$ and $\bs \mu^{\pmb *}$. Proofs of Proposition 1 and
Theorems~1,~2
\endsubhead
\medskip

We start with the

\demo{Proof of Proposition 1}  Using Corollary~B.2 in
Appendix B we see that the sequence $(u_n)$ is non-increasing.
Also (see Corollary~B.1)
$$
\| g_n(u_n)\|_{L^1} \leq \|\mu \|_{\Cal M}
$$
and thus
$$
\|\Delta u_n\|_{\Cal M} \leq 2 \|\mu \|_{\Cal M}.
$$
Consequently,
$$
\|u_n\|_{L^1} \leq C \|\mu\|_{\Cal M}.
$$
Therefore, $(u_n)$ tends in $L^1$ to a limit denoted $u^*$. By Dini's
 lemma, $g_n \uparrow g$ uniformly on compact sets; thus
$$
g_n(u_n) \to g(u^*)\quad\text{a.e. }
$$

\nd
Hence $g(u^*)\in L^1(\Omega)$, (0.11)--(0.12) hold and, by Fatou's lemma,
$$
-\int_\Omega u^* \Delta \zeta + \int_\Omega g(u^*)\zeta \leq
 \int_\Omega  \zeta \, d\mu \quad   \forall
 \zeta \in C^2_0(\overline{\Omega}), \; \zeta \geq 0\text{ in } \Omega.
$$
Therefore $u^*$ is a subsolution of (0.1).  We claim that $u^*$ is the {\it largest}\/
subsolution. Indeed let $v$ be any subsolution of (0.1).
Then
$$
-\Delta v + g_n (v) \leq - \Delta v + g(v) \leq \mu \quad \text{in } (C_0^2)^*.
$$
By comparison (see Corollary~B.2)
$$
v\leq u_n \quad \text{a.e.}
$$
and, as $n\to \infty$,
$$
v\leq u^* \quad \text{a.e.}
$$
Hence $u^*$ is the largest subsolution.
\enddemo

Recall (see [FST], or Appendix A below) that any measure $\mu$ on $\Omega$
 can be uniquely decomposed as a sum of two measures, $\mu = \mu_{\trd}
+ \mu_{\trc}$ (``d" stands for diffuse and ``c" for concentrated),
 satisfying  $|\mu_{\trd}| (A) = 0$ for every Borel set
$A\subset \Omega$ such that $\capt{(A)}=0$, and
$|\mu_{\trc}|(\Omega\setminus F) = 0$ for some Borel set $F\subset
\Omega$ such that $\capt{(F)}=0$.  Note that a measure $\mu$ is
diffuse if and only if $\mu_{\trc} = 0$, i.e., $\mu=\mu_{\trd}$.
\medskip

A key ingredient in the proof of Theorems 1 and 2 is the following
version of Kato's inequality (see [K]) due to Brezis-Ponce [BP2].

\proclaim{Theorem 7 (Kato's inequality when $\bs\Delta \bs v$ is a measure) }
Let $v\in L^1(\Omega)$ be such that $\Delta v$ is a measure on
$\Omega$.  Then, for every open set $\omega \subset\subset \Omega$,
$\Delta v^+$ is a measure on $\omega$  and the following holds:
$$
\align
(\Delta v^+)_{\trd} \geq \chi_{[v\geq 0]} (\Delta v)_{\trd} &\;
\text{ in } \omega,\tag{1.1}\\
(-\Delta v^+)_{\trc} =  (-\Delta v)^+_{\trc}\quad &\; \text{ in }  \omega.\tag{1.2}
\endalign
$$
\endproclaim

Note that the right-hand side of (1.1) is well-defined because the
function $v$ is quasi-continuous.  More precisely, if $v\in L^1(\Omega)$ and
$\Delta v$ is a measure, then there exists $\tilde v: \Omega\to \Bbb
R$ quasi-continuous such that $v=\tilde v$ a.e.$ $
in $\Omega$ (see [A1]
and also [BP1, Lemma~1]).  Recall that $\tilde v$ is quasi-continuous
if and only if, given any $\varepsilon > 0$, one can find an open set
$\omega_\eps \subset \Omega$ such that $\capt{(\omega_\eps)} <\varepsilon$ and
$\tilde v|_{\Omega\setminus \omega_\eps}$ is continuous.  In
particular, $\tilde v$ is finite
q.e. (= quasi-everywhere = outside a set of zero capacity).  It is
easy to see that $\chi_{[\tilde v\geq 0]}$ is integrable with respect
to the measure $|(\Delta v)_d|$.  When $v\in L^1$ and $\Delta v$ is a
measure, we will systematically replace $v$ by its quasi-continuous
representative.

\medskip
Here are two consequences of Theorem~7 which will be used in the
sequel.  The first one was originally established by
Dupaigne-Ponce~[DP] and it is equivalent to (1.2):

\proclaim{Corollary 8  (``Inverse'' maximum principle)}  Let $v\in
L^1(\Omega)$ be such that $\Delta v$ is a measure.
If $v\geq 0$ a.e.$ $ in $\Omega$, then
$$
(-\Delta v)_{\trc} \geq 0 \quad \text{in } \Omega.
$$
\endproclaim

Another corollary is the following

\proclaim{Corollary 9}  Let $u\in L^1(\Omega)$ be
such that $\Delta u$ is a measure.  Then,
$$
\Delta T_k(u) \leq \chi_{[u\leq k]}(\Delta u)_{\trd} + (\Delta
u)_{\trc}^+ \quad\text{in } \Cal D'(\Omega).
$$
\endproclaim

\noindent
Here, $T_k(s) = k-(k-s)^+$ for every $s \in \Bbb R$.

\demo{Proof}  Let $\omega \subset\subset \Omega$. Applying (1.1) and
(1.2) to $v=k-u$, yields
$$
(\Delta T_k (u))_{\trd} = - (\Delta v^+)_{\trd} \leq
- \chi_{[v \geq 0]}(\Delta v)_{\trd} = \chi_{[u\leq k]}(\Delta
u)_{\trd} \quad \text{in }\omega
$$
and
$$
(\Delta T_k(u))_{\trc} = (\Delta u)_{\trc}^+ \quad \text{in } \omega.
$$
Combining these two facts, we conclude that
$$
\Delta T_k(u) \leq \chi_{[u\leq k]}(\Delta u)_{\trd} + (\Delta
u)_{\trc}^+ \quad\text{in } \Cal D'(\omega).
$$
Since $\omega \subset\subset \Omega$ was arbitrary, the result follows.
\enddemo

\medskip
Let $u^*$ be the largest subsolution of (0.1), and define $\mu^* \in
\Cal M(\Omega)$ by (0.13). We have the following

\proclaim{Lemma 1} The reduced measure $\mu^*$ satisfies
$$
\mu^* \geq \mu_{\trd} - \mu_{\trc}^-.
$$
\endproclaim

\demo{Proof}
Let $(u_n)$ be the sequence constructed in Proposition 1.
By Corollary 9, we have
$$
\Delta T_k(u_n) \leq \chi_{[u_n\leq k]}(\Delta u_n)_{\trd} +  (\Delta
u_n)_{\trc}^+ \quad \text{in } \Cal D'(\Omega).\tag 1.3
$$
Since $u_n$ satisfies (0.10),
$$
(\Delta u_n)_{\trd} = g_n (u_n) - \mu_{\trd} \quad
\text{and} \quad (\Delta u_n)_{\trc} = - \mu_{\trc}.
$$
Inserting into (1.3) gives
$$
\align
-\Delta T_k (u_n)
& \geq \chi_{[u_n\leq k]} \big\{ \mu_{\trd} - g_n (u_n) \big\} -
\mu_{\trc}^-\\
& \geq \chi_{[u_n\leq k]} \mu_{\trd} - g_n(T_k(u_n)) - \mu_{\trc}^-
 \quad \text{in } \Cal  D'(\Omega).
\endalign
$$
For every $n\geq 1$ we
have $u^* \leq u_n \leq u_1$, so that
$$
[u^* \leq k] \supset [u_n \leq k] \supset [u_1 \leq k]
$$
and
$$
\chi_{[u_n\leq k]} \mu_{\trd} \geq \chi_{[u_1 \leq k]} \mu_{\trd}^+ - \chi_{[u^*\leq k]} \mu_{\trd}^-.
$$
Thus
$$
-\Delta T_k (u_n) + g_n (T_k(u_n)) \geq \chi_{[u_1 \leq k]}
 \mu_{\trd}^+ - \chi_{[u^*\leq k]} \mu_{\trd}^- - \mu_{\trc}^-
 \quad\text{in } \Cal D'(\Omega).\tag 1.4
$$
By dominated convergence,
$$
g_n(T_k(u_n)) \to g(T_k(u^*))\quad \text{in } L^1(\Omega),
\quad\text{as } n\to \infty.
$$
As $n\to \infty$ in (1.4), we get
$$
-\Delta T_k(u^*) + g(T_k(u^*))\geq \chi_{[u_1 \leq k]}
 \mu_{\trd}^+ - \chi_{[u^*\leq k]} \mu_{\trd}^- - \mu_{\trc}^- \quad
 \text{in } \Cal D' (\Omega).
$$
Let $k\to \infty$.  Since  both sets $[u_1 = + \infty]$ and  $[u^* = + \infty]$ have zero
capacity (recall that $u_1$ and $u^*$ are quasi-continuous and, in
particular, both functions are finite q.e.), we conclude that
$$
\mu^* = - \Delta u^* + g(u^*) \geq \mu_{\trd}^+ -
\mu_{\trd}^- - \mu_{\trc}^- = \mu_{\trd} - \mu_{\trc}^-.
$$
This establishes the lemma.
\enddemo

\medskip

\demo{Proof of Theorems 1 and 2}
It follows from (0.14) and Lemma~1 that
$$
\mu_{\trd} - \mu_{\trc}^- \leq \mu^* \leq \mu.
$$
By taking the diffuse parts, we have
$$
(\mu^*)_{\trd} = \mu_{\trd}.\tag 1.5
$$
Thus $\mu-\mu^*=(\mu-\mu^*)_{\trc}$, which proves Theorem 2.
\medskip

We now turn to the proof of Theorem~1. Let $\lambda$ be a
good measure $\leq \mu$. We must prove that $\lambda \leq \mu^*$. Denote by $v$
the solution of (0.1) corresponding to $\lambda$,
$$
\left\{
\alignedat2
-\Delta v + g(v) & = \lambda & & \quad \text{in } \Omega,\\
v & = 0 & & \quad \text{on } \partial \Omega.
\endalignedat
\right.
$$
By (1.5),
$$
\lambda_{\trd} \leq \mu_{\trd} = (\mu^*)_{\trd}.
$$
Since $u^*$ is the largest subsolution of (0.1), we also have
$$
v\leq u^*\quad\text{a.e.}
$$
By the ``inverse'' maximum principle,
$$
\lambda_{\trc} = (-\Delta v)_{\trc} \leq (-\Delta u^*)_{\trc} =
(\mu^*)_{\trc}.
$$
Therefore $\lambda\leq \mu^*$. This establishes Theorem~1.
\enddemo

\medskip

The following lemma will be used later on:

\proclaim{Lemma 2}  Given a measure $\mu\in \Cal M(\Omega)$, let $(u_n)$
be the sequence defined in Proposition 1.  Then,
$$
g_n(u_n)\overset{*}\to{\rightharpoonup} g(u^*) + (\mu -
\mu^*) = g(u^*) + (\mu - \mu^*)_{\trc}\quad\text{weak$^*$\/ in } \Cal M(\Omega).
$$
\endproclaim

\demo{Proof}
Let $\zeta \in C_0^2(\overline\Omega)$. For every $n \geq 1$, we have
$$
\int_\Omega g_n(u_n) \zeta = \int_\Omega
u_n \Delta \zeta + \int_\Omega \zeta \, d\mu.
$$
By Proposition~1, $u_n \to u^*$ in $L^1(\Omega)$. Thus,
$$
\lim_{n \to \infty}{\int_\Omega g_n(u_n) \zeta} =  \int_\Omega
u^* \Delta \zeta + \int_\Omega \zeta \, d\mu = \int_\Omega g(u^*)
\zeta + \int_\Omega \zeta \, d(\mu - \mu^*).
$$
In other words,
$$
g_n(u_n)\overset{*}\to{\rightharpoonup} g(u^*) + (\mu -
\mu^*) \quad\text{weak$^*$\/ in } \Cal M(\Omega).
$$
Since $(\mu^*)_\trd = \mu_\trd$, the result follows.
\enddemo

\medskip


\subhead 2. Good measures. Proofs of Theorems 4, 6
\endsubhead

\medskip
We start with

\proclaim{Lemma 3}  If $\mu$ is a good measure with solution $u$, and $u_n$ is given by $(0.10)$, then
$$
u_n \to u \quad \text{in } W_0^{1,1}(\Omega) \qquad \text{and} \qquad g_n (u_n) \to
g(u) \quad \text{in } L^1(\Omega).
$$
\endproclaim

\demo{Proof}  We have
$$
-\Delta u_n + g_n (u_n) = \mu \quad\text{and}\quad -\Delta u + g(u) =
 \mu \quad \text{in } (C_0^2)^*,
$$
so that
$$
 - \Delta(u_n - u) + g_n (u_n) - g(u) = 0  \quad \text{in } (C_0^2)^*.
$$
Thus
$$
-\Delta(u_n-u) + g_n(u_n) - g_n(u) = g(u) - g_n(u)  \quad \text{in } (C_0^2)^*.
$$
Hence, by standard estimates (see Proposition~B.3),
$$
\int_\Omega |g_n(u_n) - g_n(u) | \leq \int_\Omega |g(u) - g_n(u)|\to 0.
$$
Thus
$$
\int_\Omega |g_n(u_n) - g(u)|\leq 2 \int_\Omega |g(u) - g_n(u)|\to 0.
$$
In other words, $g_n(u_n) \to g(u)$ in $L^1(\Omega)$. This clearly
implies that $\Delta(u_n - u) \to 0$ in $L^1(\Omega)$ and thus $u_n
\to u$ in $W^{1,1}_0(\Omega)$.
\enddemo
\medskip

We now turn to the

\demo{Proof of Theorem~4}  Let $u_{1,n}, u_{2,n} \in L^1(\Omega)$ be
such that
$$
\left\{
\alignedat2
-\Delta u_{i,n} +  g_n (u_{i,n}) & = \mu_i & & \quad \text{in } \Omega,\\
u_{i,n} & = 0 &  & \quad \text{on } \partial\Omega,
\endalignedat
\right.
$$
for $i=1,2$. Since $\mu_2 \leq \mu_1$, we have
$$
u_{2,n} \leq u_{1,n} \quad \text{a.e.}
$$
Thus $g_n(u_{2,n}) \leq g_n (u_{1,n}) \to g (u_1^*)$ strongly in
$L^1$ by Lemma 3.  Hence $g_n(u_{2,n}) \to g(u_2^*)$ strongly in $L^1$
and we have
$$
-\Delta u_2^* + g(u_2^*) = \mu_2 \quad \text{in } (C_0^2)^*,
$$
i.e., $\mu_2$ is a good measure.
\enddemo

\medskip
A simple property of $\Cal G$ is

\proclaim{Proposition 2}  The set $\Cal G$ of good measures is
closed with respect to strong convergence in $\Cal M(\Omega)$.
\endproclaim

\demo{Proof}
Let $(\mu_k)$ be a sequence of good measures such that $\mu_k \to \mu$
strongly in $\Cal M(\Omega)$.  For each $k\geq 1$, let $u_k$ be such
that
$$
\left\{
\alignedat2
-\Delta u_k + g(u_k) & = \mu_k && \quad \text{in } \Omega,\\
u_k & = 0 && \quad \text{on } \partial \Omega.
\endalignedat
\right.
$$
By standard estimates (see Corollary~B.1),
$$
\int_\Omega |g(u_{k_1})-g(u_{k_2})|\leq \|\mu_{k_1} - \mu_{k_2}\|_{\Cal
M}\tag 2.1
$$
and
$$
\int_\Omega |u_{k_1}-u_{k_2}| \leq C \big\|\Delta (u_{k_1} -
u_{k_2}) \big\|_{\Cal M} \leq 2 C \|\mu_{k_1} - \mu_{k_2}\|_{\Cal
M}.\tag 2.2
$$
By (2.1) and (2.2), both $(u_k)$ and $(g(u_k))$ are Cauchy sequences
in $L^1(\Omega)$. Thus, there exist $u,v \in L^1(\Omega)$ such that
$$
u_k \to u \quad \text{and} \quad g(u_k)\to v \quad \text{in } L^1(\Omega).
$$
In particular, $v = g(u)$ a.e.  It is then easy to see
that
$$
-\Delta u + g(u) = \mu\quad \text{in } (C_0^2)^*.
$$
Thus $\mu$ is a good measure.
\enddemo
\medskip

We next present a result slightly sharper than Theorem~6:

\proclaim{Theorem 6$'$}  Let $\mu\in \Cal M(\Omega)$.  The
following conditions are equivalent:
\medskip

{\rm(a)}  $\mu$ is a good measure;
\medskip

{\rm(b)}  $\mu^+$ is a good measure;
\medskip

{\rm(c)}  $\mu_{\trc}$ is a good measure;
\medskip

{\rm(d)}  $\mu = f_0 -\Delta v_0$ in $\Cal D'(\Omega)$, for some
$f_0\in L^1$ and some $v_0\in L^1$ with $g(v_0) \in L^1$.
\medskip
\endproclaim

\demo{Proof}  (a) $\Rightarrow$ (b).  Since $\mu$ and $0$ are good
measures, it follows from Corollary~4 that $\mu^+ = \sup{\{\mu, 0\}}$ is a
good measure.

\smallskip
(b) $\Rightarrow$ (a). Since $\mu^+$ is a good measure and $\mu \leq
\mu^+$ in $\Omega$, it follows from Theorem~4 that $\mu$ is a good measure.

\smallskip
(b) $\Rightarrow$ (c). Note that we always have
$$
\mu_{\trc}\leq \mu^+. \tag 2.3
$$
Indeed, $(\mu^{+} -\mu_{\trc})_{\trd}= (\mu^+)_{\trd}\geq 0$ and
$(\mu^{+} -\mu_{\trc})_{\trc}=\mu_{\trc}^+ -\mu_{\trc} \geq 0$.

\noindent
[Here and in the sequel we use the fact that $(\mu^+)_{\trd}=(\mu_{\trd})^+$
and $(\mu^+)_{\trc}=(\mu_{\trc})^+$
which will be simply denoted $\mu^+_{\trd}$ and $\mu^+_{\trc}$].

\noindent
Since $\mu^+$ is a good measure, it follows from (2.3) and
Theorem~4 that $\mu_{\trc}$ is also a good measure.

\smallskip
(c) $\Rightarrow$ (b).  It is easy to see that, for every measure
$\lambda$,
$$
\lambda^+ = \sup{\{\lambda_{\trd},\lambda_{\trc}\}}. \tag 2.4
$$
Assume $\mu_{\trc}$ is a good measure. Since
$\mu_{\trd}$ is diffuse, Corollary~2 implies that $\mu_{\trd}$ is also a
good measure.  By Corollary~4 and (2.4), $\mu^+ = \sup{\{\mu_{\trd},
\mu_{\trc}\}}$ is a good measure as well.

\smallskip
(a) $\Rightarrow$ (d).  Trivial.

\smallskip
(d) $\Rightarrow$ (c). We split the argument into two steps.

\demo{Step 1} Proof of (d) $\Rightarrow$ (c) if $v_0$ has compact support.

\smallskip
Since $\mu = f_0 -\Delta v_0$ in $\Cal D'(\Omega)$ and $v_0$ has compact
support, we have
$$
\mu = f_0 -\Delta v_0 \quad \text{in } (C_0^2)^*.
$$
Thus, $\mu - f_0 + g(v_0)$ is a good measure. Using the equivalence (a)
$\Leftrightarrow$ (c), we conclude that $\mu_{\trc} = \big[ \mu - f_0 +
g(v_0) \big]_{\trc}$ is a good measure.
\enddemo

\demo{Step 2} Proof of (d) $\Rightarrow$ (c) completed.

\smallskip
By assumption,
$$
\mu = f_0 -\Delta v_0 \quad \text{in } \Cal D'(\Omega).
$$
In particular, we have $\Delta v_0 \in \Cal M(\Omega)$, so that $v_0
\in W^{1,p}\loc(\Omega)$, $\forall p < \frac{N}{N-1}$ (see Theorem~B.1
below). Let $(\varphi_n) \subset C_\trc^\infty(\Omega)$ be such that $0
\leq \varphi_n \leq 1$ in $\Omega$ and $\varphi_n(x) = 1$ if $d(x,
\partial\Omega) > \frac{1}{n}$. Then
$$
\varphi_n \mu = f_n - \Delta(\varphi_n v_0) \quad \text{in } \Cal D'(\Omega),
$$
where
$$
f_n = \varphi_n f_0 + 2 \nabla v_0 \cdot \nabla \varphi_n + v_0 \Delta
\varphi_n \in L^1(\Omega).
$$
Moreover, since $0 \leq g(\varphi_n v_0) \leq g(v_0)$ a.e., we have
$g(\varphi_n v_0) \in L^1(\Omega)$. Thus, by Step~1,
$$
\varphi_n \mu_\trc = (\varphi_n \mu)_\trc \in \Cal G \quad \forall n
 \geq 1.
$$
Since
$\varphi_n \mu_\trc \to \mu_\trc$ strongly in $\Cal M(\Omega)$ and $\Cal G$ is closed
with respect to the strong topology in $\Cal M(\Omega)$, we conclude
that $\mu \in \Cal G$.
\enddemo

\enddemo

\medskip

We may now strengthen Corollary~7:

\proclaim{Corollary 7$'$}
We have
$$
\Cal G + \Cal M_{\trd}(\Omega) \subset \Cal G,
$$
where $\Cal M_\trd (\Omega)$ denotes the space of diffuse measures.
\endproclaim

\demo{Proof}
Let $\mu \in \Cal G$. By Theorem~6$'$, $\mu_\trc$ is a good
measure. Thus, for any $\nu \in \Cal M_{\trd}$, $(\mu +
\nu)_\trc = \mu_\trc$ is a good measure. It follows from the
equivalence (a)
$\Leftrightarrow$ (c) in the theorem above that $\mu + \nu \in \Cal G$.
\enddemo

\proclaim{Proposition~3}  Assume
$$
g(2t) \leq C(g(t) + 1)\quad \forall t\geq 0.\tag 2.5
$$
Then the set of good measures is a convex cone.
\endproclaim

\remark{Remark 7} Assumption (2.5) is called in the literature the
$\Delta_2$-condition.  It holds if $g(t) = t^p$ for $t \geq 0$ (any $p>1$), but (2.5)
{\it fails}\/ for $g(t) = \e^t-1$.  In this case, the set of good
measures is {\it not}\/  a cone.
As we will see in Section~8, Example~5, if $N=2$, then for any $a \in \Omega$ we
have $c\delta_a \in \Cal G$ if $c > 0$ is small, but $c\delta_a \not\in
\Cal G$ if $c$ is large.
\endremark

\demo{Proof of Proposition~3}  Assume $\mu \in \Cal G$.  Clearly, it
suffices to show that $2 \mu \in \Cal G$.

Let $u$ be the solution of
$$
\left\{
\alignedat2
- \Delta u + g(u) & = \mu && \quad \text{in } \Omega,\\
u & = 0 && \quad \text{on } \partial \Omega.
\endalignedat
\right.
$$
Thus,
$$
2\mu = - \Delta(2u) + 2g(u) \quad \text{in } \Cal D'(\Omega).
$$
By (2.5), $g (2u) \in L^1$.  We can now invoke the equivalence (a)
$\Leftrightarrow$ (d) in Theorem~6$'$ to conclude that $2\mu \in \Cal G$.
\enddemo

\bigskip


\subhead 3. Some properties of the mapping $\bs\mu \pmb\mapsto
\bs\mu^{\pmb *}$\endsubhead
\medskip

We start with an easy result, which asserts that the mapping $\mu
\mapsto \mu^*$ is order preserving:

\proclaim{Proposition 4}
Let $\mu, \nu \in \Cal M(\Omega)$. If $\mu \leq \nu$, then $\mu^* \leq
\nu^*$.
\endproclaim

\demo{Proof}
Since the reduced measure $\mu^*$ is a good measure and $\mu^* \leq \mu \leq \nu$, it
follows from Theorem~1 that $\mu^* \leq \nu^*$.
\enddemo

\medskip

Next, we have

\proclaim{Theorem 8}
If $\mu_1, \mu_2 \in \Cal M(\Omega)$ are mutually singular, then
$$
(\mu_1 + \mu_2)^* = (\mu_1)^* + (\mu_2)^*. \tag 3.1
$$
\endproclaim

\demo{Proof}
Since $\mu_1$ and $\mu_2$ are mutually singular, $(\mu_1)^*$ and
$(\mu_2)^*$ are also mutually singular (by (0.16)). In particular, we have
$$
(\mu_1)^* + (\mu_2)^* \leq \big[(\mu_1)^* + (\mu_2)^* \big]^+ =
\sup{\big \{ (\mu_1)^*, (\mu_2)^* \big\}}. \tag 3.2
$$
By Corollary~4, the right-hand side of (3.2) is a good measure. It
follows from Theorem~4 that $(\mu_1)^* + (\mu_2)^*$ is also a good
measure. Since
$$
(\mu_1)^* + (\mu_2)^* \leq \mu_1 + \mu_2,
$$
we conclude from Theorem~1 that
$$
(\mu_1)^* + (\mu_2)^* \leq (\mu_1 + \mu_2)^*. \tag 3.3
$$

We now establish the reverse inequality. Assume $\lambda$ is a good
measure $\leq (\mu_1 + \mu_2)$. By Radon-Nikodym, we may decompose
$\lambda$ in terms of three measures:
$$
\lambda = \lambda_0 + \lambda_1 + \lambda_2,
$$
where $\lambda_0$ is singular with respect to $|\mu_1| + |\mu_2|$,
and, for $i=1,2$, $\lambda_i$ is absolutely continuous with respect to
$|\mu_i|$. Since $\lambda_0, \lambda_1, \lambda_2 \leq \lambda^+$,
each $\lambda_j$, $j = 0,1,2$, is a good measure. Moreover, $\lambda
\leq \mu_1 + \mu_2$ implies
$$
\lambda_0 \leq 0, \quad  \lambda_1 \leq \mu_1 \quad \text{and} \quad
\lambda_2 \leq \mu_2.
$$
Thus, in particular, $\lambda_i \leq (\mu_i)^*$ for
$i=1,2$. Therefore,
$$
\lambda = \lambda_0 + \lambda_1 + \lambda_2 \leq (\mu_1)^* +
(\mu_2)^*.
$$
Since $\lambda$ was arbitrary, we have
$$
(\mu_1 + \mu_2)^* \leq (\mu_1)^* + (\mu_2)^*. \tag 3.4
$$
Combining (3.3) and (3.4), the result follows.
\enddemo

\medskip
Here are some consequences of Theorem~8:

\proclaim{Corollary 10}
For every $\mu \in \Cal M(\Omega)$, we have
$$
(\mu^*)_\trd = (\mu_\trd)^* = \mu_\trd \quad \text{and} \quad (\mu^*)_\trc =
(\mu_\trc)^*. \tag 3.5
$$
Also,
$$
(\mu^*)^+ = (\mu^+)^* \quad \text{and} \quad (\mu^*)^- =
\mu^-. \tag 3.6
$$
\endproclaim

\demo{Proof}
Since $\mu_\trd$ is a good measure (see Corollary~2), we have $(\mu_\trd)^* =
\mu_\trd$. By Theorem~8,
$$
\mu^* = (\mu_\trd + \mu_\trc)^* = (\mu_\trd)^* + (\mu_\trc)^*.
$$
Comparison between the diffuse and concentrated parts gives
(3.5). Similarly,
$$
\mu^* = (\mu^+ - \mu^-)^* = (\mu^+)^* + (-\mu^-)^* = (\mu^+)^* -\mu^-,
$$
since every nonpositive measure is good. This identity yields (3.6).
\enddemo

More generally, the same argument shows the following:

\proclaim{Corollary 11}
Let $\mu \in \Cal M(\Omega)$. For every Borel set $E \subset \Omega$,
we have
$$
(\mu \lfloor_E)^* = \mu^* \lfloor_E. \tag 3.7
$$
\endproclaim

\noindent
Here $\mu\lfloor_E$ denotes the measure defined by $\mu\lfloor_E(A) =
\mu(A \cap E)$ for every Borel set $A \subset \Omega$.

\medskip
For simplicity, from now on we shall write $\mu_\trd^* = (\mu^*)_\trd$
and $\mu_\trc^* = (\mu^*)_\trc$.

\medskip

The following result extends Corollary~7$'$:

\proclaim{Corollary 12} For every $\mu \in \Cal M(\Omega)$ and $\nu
\in \Cal M_\trd (\Omega)$,
$$
(\mu + \nu)^* = \mu^* + \nu.
$$
\endproclaim

\demo{Proof}
By Theorem~8 and Corollary~2, we have
$$
(\mu + \nu)^* = \mu_\trc^* + (\mu_\trd + \nu)^* = \mu_\trc^* + \mu_\trd +
\nu = (\mu_\trc^* + \mu_\trd^*) + \nu = \mu^* + \nu.
$$
\enddemo

Next, we have

\proclaim{Theorem 9}  Given $\mu, \nu \in \Cal M(\Omega)$,
we have
$$
\align
\big[\inf{\{ \mu, \nu\}} \big]^* & = \inf{\{ \mu^*, \nu^*\}},\tag 3.8\\
\big[\sup{\{ \mu, \nu\}} \big]^* & = \sup{\{ \mu^*, \nu^*\}}.\tag 3.9
\endalign
$$
\endproclaim

\demo{Proof}

\noindent
{\it Step 1.} Proof of (3.8).
\smallskip

Clearly,
$$
\inf{\{ \mu^*, \nu^*\}} \leq \big[ \inf{ \{\mu, \nu\}} \big]^*.
$$
Assume $\lambda$ is a good measure $\leq \inf{\{\mu, \nu\}}$.  By
Theorem~1, $\lambda \leq \mu^*$ and $\lambda \leq \nu^*$.  Thus,
$\lambda \leq \inf{\{\mu^*, \nu^*\}}$, whence
$$
\big[ \inf{ \{\mu,\nu\}} \big]^* \leq \inf{\{\mu^*, \nu^*\}}.
$$

\noindent
{\it Step 2.} Proof of (3.9).
\smallskip

Applying the Hahn decomposition to $\mu - \nu$, we may write
$\Omega$ in terms of two disjoint Borel sets $E_1, E_2 \subset
\Omega$, $\Omega = E_1 \cup E_2$, so that
$$
\mu \geq \nu \quad \text{in $E_1$} \quad \text{and} \quad \nu \geq \mu
\quad \text{in $E_2$.}
$$
By Proposition~4 and Corollary~11,
$$
\mu^* \lfloor_{E_1} = (\mu \lfloor_{E_1})^* \geq (\nu \lfloor_{E_1})^* =
\nu^*\lfloor_{E_1}.
$$
Thus, $\mu^* \geq \nu^*$ on $E_1$. Similarly, $\nu^* \geq \mu^*$ on
$E_2$. We then have
$$
\sup{ \{\mu, \nu\} } = \mu\lfloor_{E_1} + \nu\lfloor_{E_2} \quad
\text{and} \quad \sup{ \{\mu^*, \nu^*\} } = \mu^*\lfloor_{E_1} +
\nu^*\lfloor_{E_2}. \tag 3.10
$$
On the other hand, by Theorem~8 and Corollary~11,
$$
(\mu \lfloor_{E_1} + \nu \lfloor_{E_2})^* = (\mu \lfloor_{E_1})^* +
(\nu \lfloor_{E_2})^* = \mu^* \lfloor_{E_1} + \nu^*
\lfloor_{E_2}. \tag 3.11
$$
Combining (3.10) and (3.11), we obtain (3.9).
\enddemo

\medskip

We now show that $\mu \mapsto \mu^*$ is non-expansive:

\proclaim{Theorem 10}  Given $\mu, \nu \in \Cal M (\Omega)$,
we have
$$
|\mu^* - \nu^* |\leq |\mu-\nu|. \tag 3.12
$$
More generally,
$$
(\mu^*-\nu^*)^+ \leq (\mu - \nu)^+. \tag 3.13
$$
\endproclaim

\demo{Proof}
Clearly, it suffices to show that (3.13) holds. We split
the proof into two steps.

\smallskip
\noindent
{\it Step 1.} Assume $\nu\leq \mu$. Then we claim that
$$
\mu^* - \nu^*\leq \mu - \nu. \tag 3.14
$$

Indeed, let $v_n$ be the solution of (0.10) corresponding to the measure $\nu$.
 Since $\nu \leq \mu$, we have
$$
v_n \leq u_n \quad \text{a.e.,} \quad\forall n \geq 1.
$$
Recall that $g_n$ is nondecreasing; thus,
$$
g_n(v_n) \leq g_n(u_n) \quad \text{a.e.}
$$
Let $n\to \infty$.  According to Lemma 2, we have
$$
g(v^*) + (\nu - \nu^*)_{\trc} \leq g(u^*) + (\mu - \mu^*)_{\trc}.
$$
Taking the concentrated part on both sides of this inequality yields
$$
(\nu - \nu^*)_{\trc} \leq (\mu - \mu^*)_{\trc}.
$$
Since $\nu_{\trd} = \nu^*_{\trd}$ and $\mu_{\trd} = \mu^*_{\trd}$ (by
Corollary~2), we have
$$
\nu - \nu^* \leq \mu - \mu^* ,
$$
which is (3.14).
\smallskip

\noindent
{\it Step 2.}\/  Proof of (3.1) completed.
\smallskip

Recall that
$$
\sup{\{\mu, \nu\}} = \nu + (\mu-\nu)^+.\tag 3.15
$$
Applying the previous step to the measures $\nu$ and $\sup{\{ \mu, \nu\}}$, we have
$$
\big[\sup{\{\mu, \nu\}}\big]^* - \nu^*\leq \sup{\{\mu, \nu\}} - \nu =
(\mu-\nu)^+. \tag 3.16
$$
By (3.9), (3.15) and (3.16),
$$
(\mu-\nu)^+\geq \big[\sup{ \{\mu, \nu\}} \big]^*-\nu^* = \sup{\{\mu^*,
\nu^*\}} - \nu^*= (\mu^* - \nu^*)^+.
$$
Therefore, (3.13) holds.
\enddemo

\bigskip


\subhead 4. Approximation of $\bs\mu$ by $\bs\rho_{\bs n} \pmb * \bs \mu$ \endsubhead
\medskip

Let $(\rho_n)$ be a sequence of mollifiers in $\Bbb R^N$ such
that $\supp{\rho_n} \subset B_{1/n}$ for every $n \geq 1$.
Given $\mu \in \Cal M (\Omega)$, set
$$
\mu_n= \rho_n*\mu,
$$
that is,
$$
\mu_n(x) = \int_\Omega \rho_{n}(x-y) \, d\mu(y)  \quad \forall x
\in \Bbb R^N. \tag 4.1
$$
[The integral in (4.1) is well-defined in view of Proposition~C.1 in
Appendix~C below. Here, we identify $\mu$ with $\tilde\mu \in
\big[ C(\overline\Omega) \big]^*$ defined there].

\medskip

Let $u_n$ be the solution of
$$
\left\{
\alignedat2
- \Delta u_n + g(u_n) & = \mu_n &&\quad \text{in } \Omega,\\
u_n & = 0 && \quad  \text{on } \partial \Omega.
\endalignedat
\right. \tag 4.2
$$

\proclaim{Theorem 11}  Assume in addition that $g$ is convex.
Then $u_n\to u^*$ in $L^1(\Omega)$, where $u^*$ is given by Proposition~1.
\endproclaim

\demo{Proof}

\noindent{\it Step 1.}\/ The conclusion holds if $\mu$ is a good measure.
\smallskip

In this case, there exists $u = u^*$ such that
$$
\left\{
\alignedat2
-\Delta u + g(u) & = \mu && \quad \text{in } \Omega,\\
u & = 0 && \quad \text{on } \partial \Omega.
\endalignedat
\right. \tag 4.3
$$
Let $\omega\subset\subset\Omega$.  For $n\geq 1$ sufficiently large, we have
$$
-\Delta (\rho_n * u) + \rho_n * g(u) = \mu_n\quad
\text{in } \omega.
$$
Thus, using the convexity of $g$,
$$
\Delta (\rho_n* u - u_n)= \rho_n* g(u) - g(u_n) \geq g(\rho_n * u ) -
g(u_n)\quad \text{in } \omega.
$$
By the standard version of Kato's inequality (see [K]),
$$
\Delta (\rho_n* u - u_n)^+ \geq \big\{g(\rho_n *
u) - g(u_n)\big\}^+\geq 0 \quad\text{in } \Cal D'(\omega). \tag 4.4
$$
Since
$$
\int_\Omega |\Delta u_n | \leq 2 \|\mu_n\|_{\Cal M} \leq C \quad \forall n\geq
1,
$$
we can extract a subsequence $(u_{n_k})$ such that
$$
u_{n_k} \to v \quad\text{in } L^1(\Omega),
$$
for some $v \in W_0^{1,1}(\Omega)$.
As $n_k\to \infty$ in (4.4), we have
$$
-\Delta(u-v)^+\leq 0\quad\text{in } \Cal D'(\omega).
$$
Since $\omega\subset\subset\Omega$ was arbitrary,
$$
- \Delta (u - v)^+\leq 0\quad  \text{in } \Cal D'(\Omega). \tag 4.5
$$
On the other hand,
$$
(u-v)^+ \in W^{1,1}_0(\Omega). \tag 4.6
$$
From (4.5), (4.6) and the weak form of the maximum principle (see
Proposition~B.1) we deduce that
$$
(u - v)^+ \leq 0 \quad\text{a.e.}
$$
Therefore,
$$
v\geq u \quad \text{a.e.}
$$
By Fatou's lemma, $v$ is a subsolution of (0.1); comparison with (4.3) yields,
$$
v\leq u \quad \text{a.e.}
$$
We conclude that
$$
v = u \quad \text{a.e.}
$$
Since $v$ is independent of the subsequence $(u_{n_k})$, we must have
$$
u_n\to u = u^* \quad \text{in } L^1(\Omega).
$$

\noindent{\it Step 2.}\/  Proof of Theorem 11 completed.
\smallskip

Without loss of generality, we may assume that
$$
u_n \to v \quad\text{in } L^1(\Omega).
$$
By Fatou, once more, $v$ is a subsolution of (0.1). Proposition 1 yields
$$
v\leq u^* \quad \text{a.e.}
$$
Let $u^*_n$ denote the solution of
$$
\left\{
\alignedat2
- \Delta u_n^* + g (u^*_n) & = \rho_{n}* \mu^* && \quad \text{in
} \Omega,\\
u^*_n & = 0  && \quad \text{on } \partial \Omega.
\endalignedat
\right.
$$
By the previous step,
$$
u^*_n \to u^* \quad\text{in } L^1(\Omega).
$$
On the other hand, we know from the maximum principle that
$$
u_n^* \leq u_n \quad \text{a.e.}
$$
Thus, as $n\to \infty,$
$$
u^*\leq v \quad \text{a.e.}
$$
Since $v\leq u^*$ a.e., the result follows.
\enddemo

\medskip

\definition{Open problem 1}
Does the conclusion of Theorem 11 remain valid without the convexity
assumption on $g$~?
\enddefinition

\bigskip


\subhead 5.  Further convergence results\endsubhead
\medskip

We start with the following

\proclaim{Theorem 12}  Let $(f_n) \subset L^1(\Omega)$ and $f\in
L^1(\Omega)$.  Assume
$$
f_n \rightharpoonup f\quad \text{weakly in } L^1. \tag 5.1
$$
Let $u_n$  (resp. $u$) be the solution of\/ $(0.1)$ associated with $f_n$
(resp. $f$).  Then $u_n \to u$ in $L^1(\Omega)$.
\endproclaim

\demo{Proof}
By definition,
$$
-\Delta u_n + g(u_n) = f_n \quad \text{and} \quad -\Delta u + g(u) = f
 \quad \text{in } (C_0^2)^*.
$$
Using a device introduced by Gallou\"et-Morel~[GM] (see also
Proposition~B.2 below), we
have, for every $M>0$,
$$
\underset{[|u_n|\geq M]}\to\int |g(u_n)| \leq \underset{[|u_n|\geq M]}\to\int |f_n|.
$$
Thus
$$
\int_E |g(u_n)|= \underset{[|u_n|\geq M]}\to{\int_E} + \underset{[|u_n|<
M]}\to{\int_E}\leq \underset{[|u_n|\geq M]}\to{\int} |f_n| + g(M)|E|. \tag 5.2
$$
On the other hand, $\|\Delta u_n \|_{L^1}  \leq C $ implies
$\|u_n\|_{L^1} \leq C$, and thus
$$
\meas{\big[|u_n|\geq M \big]} \leq \frac{C}{M}.
$$
From (5.1) and a theorem of Dunford-Pettis (see, e.g., [DS,
Corollary~IV.8.11]) we infer that $(f_n)$ is equi-integrable. Given
$\delta >
0$,  {\it fix} $M > 0$ such that
$$
\underset{[|u_n|\geq M]}\to{\int} |f_n| \leq \delta \quad \forall n
\geq 1. \tag 5.3
$$
With this fixed $M$, choose $|E|$ so small that
$$
g(M) |E| < \delta. \tag 5.4
$$
We deduce from (5.2)--(5.4) that $g(u_n)$ is equi-integrable.

\nd
Passing to a subsequence, we may assume that $u_{n_k} \to v$ in
$L^1(\Omega)$ and a.e., for some $v\in L^1(\Omega)$. Then $g(u_{n_k})
\to g(v)$ a.e. By Egorov's lemma, $g(u_{n_k}) \to g(v)$ in
$L^1(\Omega)$. It follows that $v$ is a solution of (0.1) associated
to $f$. By the uniqueness of the limit, we must have $u_n \to u$ in
$L^1(\Omega)$.
\enddemo

\medskip
\remark{Remark 8} Theorem~12 is no longer true if one replaces the
weak convergence $f_n
\rightharpoonup f$ in $L^1$, by the weak$^*$ convergence in the sense
of measures. Here is an example:

\medskip
\example{Example 1}
Assume $N \geq 3$ and let $g(t)=(t^{+})^{q}$ with $q \geq
\frac{N}{N-2}$. Let $f\equiv1$ in $\Omega$. We will construct a
sequence $(f_k)$ in $C_\trc^\infty(\Omega)$ such that
$$
f_k \overset{*}\to\rightharpoonup f \quad \text{in } \Cal M(\Omega), \tag 5.5
$$
and such that the solutions $u_k$ of (0.1) corresponding to $f_k$
converge to $0$ in $L^1(\Omega)$.
Let $(\mu_k)$ be any sequence in $\Cal M(\Omega)$ converging weak$^*$
to $f$, as   $k\to \infty$, and such that each measure $\mu_k$ is a
linear combination of
Dirac masses. (For example, each $\mu_k$ can be of the form
$|\Omega|M^{-1}\sum \delta_{a_i}$, where the $M$ points $a_i$ are
uniformly distributed in $\Omega$). Recall that for $ \mu=\delta_a$,
the corresponding $u^*$ in Proposition 1 is $\equiv0$ (see [B4]
or Theorem~B.6 below). Similarly, for each $\mu_k$, the
corresponding $u^*$ is $\equiv0$. Set
$h_{n,k}=\rho_{n}*\mu_k$, with the same notation as in
Section~4. Let $u_{n,k}$ denote the solution of (0.1) relative to
$h_{n,k}$. For each {\it fixed} $k$ we know, by Theorem~11, that
$u_{n,k} \to 0$ strongly in $L^1(\Omega)$ as   $n\to \infty$. For each
$k$, choose $N_{k}>k$ sufficiently large so that
$\|u_{N_{k},k}\|_{L^1}<1/k$. Set $f_k=h_{N_{k},k}$, so that
$u_k=u_{N_{k},k}$ is the corresponding solution of (0.1). It is easy
to check that, as $k\to \infty$,
$$
f_k \overset{*}\to\rightharpoonup f \equiv1 \quad \text{in } \Cal M(\Omega),
\quad \text{but} \quad u_k \to 0 \quad \text{in } L^1(\Omega).
$$
\endexample
\endremark

Our next result is a refinement of Theorem~12 in the spirit of
Theorem~6. Let  $\mu \in \Cal M(\Omega)$ and let $(\mu_n)$ be a
sequence in $\Cal M(\Omega)$. Assume that
$$
\alignat2
\mu & = f - \Delta v &&  \quad \text{in } (C_0^2)^*, \tag 5.6\\
\mu_n & = f_n - \Delta v_n &&  \quad \text{in } (C_0^2)^*, \tag 5.7
\endalignat
$$
where $f \in L^1$, $f_n \in L^1$, $v\in L^1$, $v_n \in L^1$, $g(v)\in
L^1$, and $g(v_n)\in L^1$.

\medskip
By Theorem~6 we know that there exist $u$ and $u_n$  solutions of
$$
\alignat4
-\Delta u + g(u) & = \mu && \quad \text{in } \Omega, & \quad u & = 0
&& \quad \text{on } \partial \Omega, \tag 5.8\\
-\Delta u_n + g(u_n) & = \mu_n &&\quad \text{in } \Omega,& \quad u_n &
= 0 && \quad \text{on } \partial \Omega. \tag 5.9
\endalignat
$$

\proclaim{Theorem~13}  Assume $(5.6)$--$(5.9)$ and moreover
$$
\align
&\|\mu_n \|_{\Cal M}\leq C , \tag 5.10\\
&f_n \rightharpoonup f \quad \text{weakly in } L^1, \tag 5.11\\
&v_n \to v \quad \text{in } L^1 \quad \text{and} \quad g(v_n) \to g(v)
\quad \text{in } L^1. \tag 5.12
\endalign
$$
Then $u_n \to u$ in $L^1(\Omega)$.
\endproclaim

\demo{Proof} We divide the proof into two steps.
\smallskip

\noindent
{\it Step 1.} Fix $0 < \alpha < 1$ and let $u(\alpha)$, $u_n(\alpha)$ be
the solutions of
$$
\alignat4
-\Delta  u(\alpha) + g(u(\alpha)) &= \alpha \mu \quad &&\text{in }
 \Omega, \quad  & u(\alpha) & = 0
\quad &\text{on } \partial \Omega, \tag 5.13\\
-\Delta  u_n(\alpha) + g(u_n(\alpha)) &= \alpha \mu_n \quad &&\text{in } \Omega,
 \quad  &u_n(\alpha) &= 0
\quad &\text{on } \partial \Omega .\tag 5.14
\endalignat
$$
Then $u_n(\alpha) \to  u(\alpha)$ in $L^1(\Omega)$.

\medskip
Note that  $u(\alpha)$ and  $u_n(\alpha)$ exist since $\alpha\mu =
\alpha f - \Delta(\alpha v)$ and $g(\alpha v)\leq g(v)$, so
that  $g(\alpha v)\in L^1$, and similarly for $\alpha\mu_n$. We may
then apply Theorem~6 once more. For simplicity we will omit the
dependence in $\alpha$ and we will write $\tilde u$, $\tilde u_n$
instead of $ u(\alpha)$, $u_n(\alpha)$ (recall that in this step
$\alpha$ is {\it fixed}\/).
Since
$$
\|\Delta\tilde u_n\|_{\Cal M} \leq 2\alpha\|\mu_n\|_{\Cal M} \leq C,
$$
we can extract a subsequence of $(\tilde u_n)$ converging strongly in
$L^1(\Omega)$ and a.e. Let $w \in W^{1,1}_0(\Omega)$ be such
that $\tilde u_{n_k} \to w$ in $L^1(\Omega)$ and a.e. We will
prove that $w$ satisfies (5.13), and therefore, by uniqueness, $w
= \tilde u$. Since
$w$ is independent of the subsequence, we will infer that $(\tilde
u_n)$ converges to $\tilde u$, which is the desired conclusion.

\medskip
We claim that
$$
g(\tilde u_n) \quad \text{is equi-integrable}. \tag5.15
$$
To establish (5.15) we argue as in the proof of Theorem~12. From
(5.7) and (5.14) we see that
$$
-\Delta(\tilde u_n - \alpha v_n) + [g(\tilde u_n) - g(\alpha v_n)]=
 h_n \quad \text{in } (C_0^2)^*, \tag5.16
$$
with
$$
h_n=\alpha f_n - g(\alpha v_n). \tag 5.17
$$
Using (5.11) and (5.12) we see that
$$
(h_n)  \quad \text{is equi-integrable}. \tag5.18
$$
From (5.16) and Proposition~B.2 we obtain (as in the proof of Theorem~12)
that, for every $M>0$,
$$
\underset{[|\tilde u_n-\alpha v_n|\geq M]}\to\int |g(\tilde u_n)-g(\alpha v_n)| \leq \underset{[|\tilde u_n-\alpha v_n|\geq M]}\to\int |h_n|. \tag 5.19
$$
On the other hand, for any Borel set $E$ of $\Omega$, we have
$$
\int_E g(\tilde u_n) = \int_{A_n} g(\tilde u_n) + \int_{B_n} g(\tilde
u_n) + \int_{C_n} g(\tilde u_n), \tag 5.20
$$
where
$$
\align
A_n & = [\tilde u_n \geq v_n] \cap [|\tilde u_n-\alpha v_n|\geq M] \cap E,\\
B_n & = [\tilde u_n \geq v_n] \cap [|\tilde u_n-\alpha v_n|< M] \cap E,\\
C_n & = [\tilde u_n < v_n] \cap E.
\endalign
$$
To handle the integral on $A_n$, write
$$
\int_{A_n} g(\tilde u_n) \leq \underset{[|\tilde u_n-\alpha v_n|\geq M]}\to\int |g(\tilde u_n)-g(\alpha v_n)| + \int_E g(v_n).
$$
Thus, by (5.19),
$$
\int_{A_n} g(\tilde u_n) \leq \underset{[|\tilde u_n-\alpha v_n|\geq
M]}\to\int |h_n| + \int_E g(v_n) .\tag 5.21
$$
Next, on $B_n$, we have
$$
\tilde u_n< M +\alpha v_n \leq M + \alpha \tilde u_n
$$
and thus
$$
\tilde u_n < \frac{M}{1-\alpha}.
$$
Therefore
$$
\int_{B_n} g(\tilde u_n) \leq g\Big( \frac{M}{1-\alpha} \Big) |E|. \tag 5.22
$$
Finally we have
$$
\int_{C_n} g(\tilde u_n) \leq \int_E g(v_n) . \tag 5.23
$$
Combining (5.20)--(5.23) yields
$$
\int_E g(\tilde u_n) \leq  \underset{[|\tilde u_n-\alpha v_n|\geq
M]}\to\int |h_n| + 2\int_E g(v_n) + g \Big( \frac{M}{1-\alpha} \Big) |E| .\tag 5.24
$$
But $\|\tilde u_n-\alpha v_n\|_{L^1} \leq C $ and therefore
$$
\meas{[|\tilde u_n-\alpha v_n|\geq M]} \leq \frac{C}{M}. \tag 5.25
$$
Given $\delta > 0$,  {\it fix} $M > 0$ sufficiently large such that
$$
 \underset{[|\tilde u_n-\alpha v_n|\geq M]}\to\int |h_n| \leq \delta
 \quad \forall n\geq1
$$
(here we use (5.18) and (5.25)).
With this fixed $M$, choose $|E|$ so small that
$$
2\int_E g(v_n) + g \Big( \frac{M}{1-\alpha} \Big) |E| \leq \delta
\quad \forall n\geq1.
$$
This finishes the proof of (5.15).

Since $g(\tilde u_n) \to g(w)$ a.e., we deduce from (5.15) and Egorov's lemma that
$g(\tilde u_n) \to g(w)$ in $L^1$.  We are now able to pass to the
limit in (5.14) and conclude that $w$ satisfies (5.13), which was the
goal of Step 1.

\medskip
\noindent
{\it Step 2.} Proof of the theorem completed.
\smallskip

Here the dependence on $\alpha$ is important and we return to the notation $u(\alpha)$ and $u_n(\alpha)$.
From (5.8) and (5.13) we deduce that
$$
\|\Delta(u(\alpha)-u)\|_{\Cal M} \leq 2(1-\alpha)\|\mu\|_{\Cal M} \tag 5.26
$$
and similarly, from (5.9) and (5.14), we have
$$
\|\Delta(u_n(\alpha)-u_n)\|_{\Cal M} \leq 2(1-\alpha)\|\mu_n\|_{\Cal
M} \leq C(1-\alpha). \tag 5.27
$$
Estimates (5.26) and (5.27) yield
$$
\|u(\alpha)-u\|_{L^1} + \|u_n(\alpha)-u_n\|_{L^1} \leq C(1-\alpha), \tag 5.28
$$
with $C$ independent of $n$ and $\alpha$.
Finally we write
$$
\|u_n - u\|_{L^1} \leq \|u(\alpha)-u\|_{L^1} + \|u_n(\alpha)-u_n\|_{L^1} + \|u_n(\alpha)-u(\alpha)\|_{L^1}. \tag 5.29
$$
Given $\varepsilon>0$, fix $\alpha<1$ so small that
$$
C(1-\alpha)<\varepsilon \tag 5.30
$$
and then apply Step 1 to assert that
$$
\|u_n(\alpha)-u(\alpha)\|_{L^1} < \varepsilon \quad \forall n \geq N, \tag5.31
$$
provided $N$ is sufficiently large.
Combining (5.28)--(5.31) yields
$$
\|u_n - u\|_{L^1} \leq 2\varepsilon \quad\forall n \geq N,
$$
which is the desired conclusion.
\enddemo

\bigskip


\subhead 6.  Nonnegative measures which are good for every $\bs g$ must be diffuse
\endsubhead
\medskip

Let $h : [0,\infty) \to [0,\infty)$ be a continuous nondecreasing
function with $h(0) = 0$. Given a compact set $K \subset \Omega$, let
$$
\capt_{\Delta, h}(K) = \inf{\left\{ \int_\Omega h(|\Delta \varphi|) \;
; \; \varphi \in C_\trc^\infty(\Omega), \; 0 \leq \varphi \leq 1, \;
  \text{and } \varphi = 1 \text{ on } K \right\}},
$$
where, as usual, $C_\trc^\infty(\Omega)$ denotes the set of
$C^\infty$-functions with compact support in $\Omega$.

We start with

\proclaim{Proposition 5}
Assume
$$
\lim_{t \to \infty}{\frac{g(t)}{t}} = + \infty \quad \text{and} \quad
g^*(s) > 0 \quad \text{for } s >0 . \tag 6.1
$$
If $\mu$ is a good measure, then $\mu^+(K) = 0$ for every compact set
$K \subset \Omega$ such that $\capt_{\Delta,g^*}(K) = 0$.
\endproclaim

Here, $g^*$ denotes the convex conjugate of $g$, which is
finite in view of the coercivity of $g$. Note that if $g'(0) =
0$, then $g^*(s) > 0$ for every $s > 0$.
\medskip

\demo{Proof}
Since $\mu$ is a good measure, $\mu^+$ is also a good measure. Thus,
$$
\mu^+ = - \Delta v + g(v) \quad \text{in } (C_0^2)^*
$$
for some $v \in L^1(\Omega)$, $v \geq 0$ a.e., such that $g(v) \in L^1(\Omega)$.

\nd
Let $\varphi_n \in C_\trc^\infty(\Omega)$ be such that $0 \leq
\varphi_n \leq 1$ in $\Omega$,
$\varphi_n = 1$ on $K$, and
$$
\int_\Omega g^*(|\Delta \varphi_n|) \to 0.
$$
Passing to a subsequence if necessary, we may assume that
$$
g^*(|\Delta \varphi_n|) \to 0 \quad \text{a.e.} \quad
\text{and} \quad g^*(|\Delta \varphi_n|) \leq G \in L^1(\Omega) \quad
\forall n \geq 1.
$$
Since $g^*(s) > 0$ if $s > 0$, we also have
$$
\varphi_n,\, |\Delta \varphi_n| \to 0 \quad \text{a.e.}
$$

\nd
For every $n \geq 1$,
$$
\mu^+(K) \leq \int_\Omega \varphi_n \, d\mu^+ =
\int_\Omega \big[ g(v) \varphi_n - v \Delta \varphi_n \big]. \tag 6.2
$$
Note that
$$
\big| g(v) \varphi_n - v \Delta \varphi_n \big| \to 0 \quad \text{a.e.}
$$
and
$$
\big| g(v) \varphi_n - v \Delta \varphi_n \big| \leq  2 g(v) + g^*(|\Delta
 \varphi_n|) \leq 2 g(v) + G \in L^1(\Omega).
$$
By dominated convergence, the right-hand side of (6.2) converges to 0
as $n \to \infty$. We then conclude that $\mu^+(K) = 0$.
\enddemo
\medskip

As a consequence of Proposition~5 we have

\proclaim{Theorem~14}
Given a Borel set $\Sigma \subset \Omega$ with zero $H^1$-capacity,
there exists $g$ such that
$$
\mu^* = - \mu^- \quad \text{for every measure $\mu$ concentrated on
$\Sigma$.}
$$
In particular, for every nonnegative $\mu \in \Cal M(\Omega)$
concentrated on a set of zero $H^1$-capacity, there exists some $g$ such
that $\mu^* = 0$.
\endproclaim

\demo{Proof}
Let $\Sigma \subset \Omega$ be a Borel set of zero $H^1$-capacity.
Let $(K_n)$ be an increasing sequence of compact sets in $\Sigma$ such
that
$$
\mu^+ \Big( \Sigma \setminus \bigcup_n{K_n} \Big) = 0.
$$
For each $n \geq 1$, $K_n$ has zero $H^1$-capacity. By Lemma~E.1, one can find $\psi_n \in
C_\trc^\infty(\Omega)$ such that $0 \leq \psi_n \leq 1$ in $\Omega$, $\psi_n = 1$
in some neighborhood of $K_n$, and
$$
\int_\Omega |\Delta \psi_n| \leq \frac{1}{n} \quad \forall n \geq
1.
$$
In particular, $\Delta \psi_n \to 0$ in $L^1(\Omega)$.
Passing to a subsequence if necessary, we may assume that
$$
\Delta \psi_n \to 0 \quad \text{a.e.} \quad \text{and} \quad |\Delta \psi_n| \leq G \in L^1(\Omega) \quad \forall n \geq 1.
$$
According to a theorem of De La Vall\'ee-Poussin (see [DVP,
Remarque~23] or [DM,
Th\'eor\`eme~II.22]), there exists a convex
function $h : [0,\infty) \to [0, \infty)$ such that $h(0) = 0$, $h(s)
> 0$ for $s > 0$,
$$
\lim_{t \to \infty}{\frac{h(t)}{t}} = + \infty, \quad \text{and} \quad
h(G) \in L^1(\Omega).
$$
By dominated convergence, we then have $h(|\Delta \psi_n|) \to 0$ in $L^1(\Omega)$.
Thus,
$$
\capt_{\Delta, h} (K_n) = 0 \quad \forall n \geq 1. \tag 6.3
$$
Let $g(t) = h^*(t)$ if $t \geq 0$, and $g(t) = 0$ if $t < 0$. By
duality, $h = g^*$ on $[0, \infty)$.

\nd
Let $\mu \in \Cal M(\Omega)$ be any measure concentrated on $\Sigma$. By
Proposition~5, (6.3) yields
$$
(\mu^*)^+(K_n) = 0 \quad \forall n \geq 1,
$$
where the reduced measure $\mu^*$ is computed with respect to $g$.
Thus, $(\mu^*)^+(\Sigma) = 0$. Since $\mu$ is concentrated on
$\Sigma$, we have $(\mu^*)^+ = 0$. Applying Corollary~10, we then get
$$
\mu^* = (\mu^*)^+ - (\mu^*)^- = - \mu^-,
$$
which is the desired result.
\enddemo

We may now present the

\demo{Proof of Theorem~5}
Assume $\mu \in \Cal M(\Omega)$ is a good
measure for every $g$. Given a Borel set $\Sigma \subset \Omega$ with
zero $H^1$-capacity, let $\lambda = \mu^+ \lfloor_\Sigma$. In view of
Theorem~14, there exists $\tilde g$ for which $\lambda^* = 0$. On the
other hand, by
Theorems~4 and 6$'$, $\lambda$ is a good measure for $\tilde g$. Thus,
$\lambda = \lambda^* = 0$. In other words, $\mu^+(\Sigma) = 0$. Since
$\Sigma$ was arbitrary, $\mu^+$ is diffuse.
This establishes the theorem.
\enddemo

\smallskip
We conclude this section with the following

\definition{Open problem 2}
Let $g : \Bbb R \to \Bbb R$ be any given continuous,
nondecreasing function satisfying (0.5).
Can one always find some nonnegative $\mu \in \Cal M(\Omega)$ such that $\mu$ is
good for $g$, but $\mu$ is {\it not}\/ diffuse?
\enddefinition

After this paper was finished, A.C.~Ponce~[P] has given a positive
answer to the above problem.

\bigskip


\subhead 7.  Signed measures and general nonlinearities $\bs g$
\endsubhead
\medskip

Suppose that $g : \Bbb R \to \Bbb R$ is a continuous, nondecreasing
function, such that $g(0) = 0$. But we will {\it not}\/ impose in this
section that $g(t) = 0$ if $t < 0$. We shall follow the same
approximation scheme as in the Introduction. Namely, let $(g_n)$ be a
sequence of nondecreasing
continuous functions, $g_n : \Bbb R \to \Bbb R$, $g_n(0) = 0$,
satisfying (0.8), such that both $(g_n^+)$ and $(g_n^-)$ verify (0.7),
and
$$
g_n^+(t) \uparrow g^+(t), \quad g_n^-(t) \uparrow g^-(t) \quad \forall
t \in \Bbb R \qquad \text{as } n \uparrow \infty.
$$

Let $\mu \in \Cal M(\Omega)$. For each $n \geq 1$, we denote by $u_n$
the unique solution of
$$
\left\{
\alignedat2
- \Delta u_n + g_n(u_n) & = \mu   && \quad \text{in } \Omega,\\
u_n & = 0  && \quad \text{on } \partial\Omega.
\endalignedat
\right. \tag 7.1
$$

First a simple observation:

\proclaim{Lemma 4}
Assume $\mu \geq 0$ or $\mu \leq 0$. Then there exists $u^* \in
L^1(\Omega)$ such that $u_n \to u^*$ in $L^1(\Omega)$. If $\mu \geq
0$, then $u^* \geq 0$ is the largest subsolution of\/ \rom{(0.1)}. If $\mu
\leq 0$, then $u^* \leq 0$ is the smallest supersolution of \/
\rom{(0.1)}. In both cases, we have
$$
\Big|\int_\Omega u^*\Delta \zeta\Big| \leq 2 \|\mu\|_{\Cal
M}\|\zeta\|_{L^\infty}\quad \forall \zeta \in C^2_0
(\overline\Omega) \tag 7.2
$$
and
$$
\int_\Omega |g(u^*)| \leq \|\mu\|_{\Cal M}. \tag 7.3
$$
\endproclaim

\demo{Proof}
If $\mu \geq 0$, then $u_n \geq 0$ a.e. In particular,
$g_n(u_n) = g_n^+(u_n)$ for every $n \geq 1$. Since $(g_n^+)$
satisfies the assumptions of Proposition~1, we conclude that $u_n \to
u^*$ in $L^1(\Omega)$, where $u^* \geq 0$ is the largest subsolution
of (0.1).

\nd
If $\mu \leq 0$, then $u_n \leq 0$, so that $w_n = -u_n$ satisfies
$$
\left\{
\alignedat2
- \Delta w_n + \tilde g_n(w_n) & = - \mu   && \quad \text{in } \Omega,\\
w_n & = 0         && \quad \text{on } \partial\Omega,
\endalignedat
\right.
$$
where $\tilde g_n(t) = g_n^-(-t)$, $\forall t \in \Bbb R$. Clearly,
the sequence $(\tilde g_n)$ satisfies the
assumptions of Proposition~1. Therefore, $u_n = - w_n \to - w^* = u^*$
in $L^1(\Omega)$. It is easy to see that $u^* \leq 0$ is the smallest
supersolution of (0.1).
\enddemo

Given $\mu \in \Cal M(\Omega)$ such that $\mu \geq 0$ or $\mu \leq 0$,
we {\it define}\/ $\mu^* \in \Cal M(\Omega)$ by
$$
\mu^* = - \Delta u^* + g(u^*) \quad \text{in } (C_0^2)^*. \tag 7.4
$$
The reduced measure $\mu^*$ is well-defined because of (7.2) and (7.3).
It is easy to see that

\medskip
(a) if $\mu \geq 0$, then $0 \leq \mu^* \leq \mu$;

\medskip
(b) if $\mu \leq 0$, then $\mu \leq \mu^* \leq 0$.

\medskip
We now consider the general case of a signed measure $\mu \in \Cal M(\Omega)$.
In view of (7.4), both measures $(\mu^+)^*$ and $(-\mu^-)^*$ are
well-defined. Moreover,
$$
-\mu^- \leq (-\mu^-)^* \leq 0 \leq (\mu^+)^* \leq \mu^+.
$$

The convergence of the approximating sequence $(u_n)$ is governed by the following:

\proclaim{Theorem 15}
Let $u_n$ be given by \rom{(7.1)}. Then, $u_n \to u^*$ in $L^1(\Omega)$,
where $u^*$ is the unique solution of
$$
\left\{
\alignedat2
- \Delta u^* + g(u^*) & = (\mu^+)^* + (-\mu^-)^*  && \quad \text{in } \Omega,\\
u^* & = 0                 && \quad \text{on } \partial\Omega.
\endalignedat
\right. \tag 7.5
$$
\endproclaim

\demo{Proof}
By standard estimates, $\| \Delta u_n \|_{\Cal M} \leq 2 \|\mu\|_{\Cal
M}$. Thus, without loss of generality, we may assume that for a
subsequence, still denoted $(u_n)$, $u_n \to u$
in $L^1(\Omega)$ and a.e. We shall show that $u$ satisfies (7.5); by
uniqueness (see Corollary~B.1), this will imply that $u$ is
independent of the subsequence.

\nd
For each $n \geq 1$, let $v_n, \tilde v_n$ be the solutions of
$$
\left\{
\alignedat2
- \Delta v_n + g_n(v_n) & = \mu^+ && \quad \text{in } \Omega,\\
v_n  & = 0  && \quad \text{on } \partial\Omega,
\endalignedat
\right. \tag 7.6
$$
and
$$
\left\{
\alignedat2
- \Delta \tilde v_n + g_n^+(\tilde v_n) & = \mu  && \quad \text{in } \Omega,\\
\tilde v_n & = 0              && \quad \text{on } \partial\Omega,
\endalignedat
\right. \tag 7.7
$$
so that $v_n \geq 0$ a.e., $v_n \downarrow v^*$ and $\tilde v_n \downarrow \tilde
v^*$ in $L^1(\Omega)$. By comparison (see Corollary~B.2), we have
$$
\tilde v_n \leq u_n \leq v_n \quad \text{a.e.}
$$
Thus,
$$
g_n^+(\tilde v_n) \leq g_n^+(u_n) \leq g_n^+(v_n) =  g_n(v_n) \quad
\text{a.e.}\tag 7.8
$$
By Lemma~2, we know that
$$
\align
g_n^+(\tilde v_n) & \overset{*}\to\rightharpoonup g^+(\tilde v^*) + \mu
- \mu^*, \\
g_n(v_n) & \overset{*}\to\rightharpoonup g(v^*) + \mu^+
- (\mu^+)^*.
\endalign
$$
Here, both reduced measures $\mu^*$ and $(\mu^+)^*$ are computed with
respect to the nonlinearity $g^+$; in particular, (see Corollary~10)
$$
\mu - \mu^* = \mu^+ - (\mu^+)^*.\tag 7.9
$$
We claim that
$$
g_n^+(u_n) \overset{*}\to\rightharpoonup g^+(u) + \mu^+ - (\mu^+)^*. \tag 7.10
$$
This will be a consequence of the following

\proclaim{Lemma 5}
Let $a_n, b_n, c_n \in L^1(\Omega)$ be such that
$$
a_n \leq b_n \leq c_n \quad \text{a.e.}
$$
Assume that $a_n \to a$, $b_n \to b$ and $c_n \to c$ a.e.$ $ in
$\Omega$ for some $a,b,c \in L^1(\Omega)$. If $(c_n - a_n)
\overset{*}\to\rightharpoonup (c - a)$
weak$^*$ in $\Cal M(\Omega)$, then
$$
(c_n - b_n) \overset{*}\to\rightharpoonup (c - b) \quad \text{weak$^*$
in $\Cal M(\Omega)$.} \tag 7.11
$$
\endproclaim

\demo{Proof}
Since
$$
0 \leq (c_n - b_n) \leq (c_n - a_n) \quad \text{a.e.,} \tag 7.12
$$
the sequence $(c_n - b_n)$ is bounded in $L^1(\Omega)$. Passing to a subsequence if
necessary, we may assume that there exists $\lambda \in \Cal
M(\Omega)$ such that
$$
(c_n - b_n) \overset{*}\to\rightharpoonup \lambda.
$$
By (7.12), we have $0 \leq \lambda \leq (c-a)$. Thus, $\lambda$ is
absolutely continuous with respect to the Lebesgue measure. In other
words, $\lambda \in L^1(\Omega)$. Given $M > 0$, we denote by $S_M$
the truncation operator $S_M(t) = \min{\{ M,
\max{\{t, -M \}} \}}$, $\forall t \in \Bbb R$. By dominated convergence, we have
$$
S_M(a_n) \to S_M(a)  \quad \text{strongly in } L^1(\Omega),
$$
and similarly for $S_M(b_n)$ and $S_M(c_n)$.
Since
$$
0 \leq \big[(c_n - S_M(c_n)) - (b_n - S_M(b_n)) \big] \leq  \big[(c_n
- S_M(c_n)) - (a_n - S_M(a_n)) \big] \quad \text{a.e.,}
$$
as $n \to \infty$ we get
$$
0 \leq \lambda - (S_M(c) - S_M(b)) \leq  \big[(c - S_M(c)) - (a -
S_M(a)) \big] \quad \text{a.e.}
$$
Let $M \to \infty$ in the expression above. We then get $\lambda
= (c - b)$. This concludes the proof of the lemma.
\enddemo

\medskip
We now apply the previous lemma with $a_n = g_n^+(\tilde v_n)$, $b_n =
g_n^+(u_n)$ and $c_n = g_n(v_n)$. In view of (7.8) and (7.9), the
assumptions of Lemma~5 are satisfied. It follows from (7.11) that
$$
g_n(v_n) - g_n^+(u_n) \overset{*}\to\rightharpoonup g(v^*) - g^+(u). \tag 7.13
$$
Thus,
$$
g_n^+(u_n) = g_n(v_n) - \big[ g_n(v_n) - g_n^+(u_n) \big]
\overset{*}\to\rightharpoonup  g^+(u) + \mu^+ - (\mu^+)^*,
$$
which is precisely (7.10). A similar argument shows that
$$
g_n^-(u_n) \overset{*}\to\rightharpoonup g^-(u) + \mu^- +
(-\mu^-)^*. \tag 7.14
$$
We conclude from (7.10) and (7.14) that
$$
g_n(u_n) \overset{*}\to\rightharpoonup g(u) + \mu - \big[ (\mu^+)^* +
(-\mu^-)^* \big]. \tag 7.15
$$
Therefore, $u$ satisfies (7.5), so that (7.5) has a solution $u^*
= u$. By uniqueness, the whole sequence $(u_n)$ converges to $u^*$ in $L^1(\Omega)$.
\enddemo

\medskip
Motivated by Theorem~15, for any $\mu \in \Cal M(\Omega)$, we {\it
define}\/ the reduced measure $\mu^*$ by
$$
\mu^* = (\mu^+)^* + (- \mu^-)^*. \tag 7.16
$$
[This definition is coherent if $\mu$ is either a positive or a
negative measure].

One can derive a number of properties satisfied by $\mu^*$. For
instance, the statements of Theorems~8--10 remain true. Moreover,

\proclaim{Theorem 2$'$}
There exists a Borel set $\Sigma \subset \Omega$ with $\capt{(\Sigma)}
= 0$ such that
$$
|\mu - \mu^*|(\Omega \backslash \Sigma) = 0.
$$
\endproclaim

\smallskip


\subhead 8. Examples\endsubhead
\medskip

We describe here some simple examples, where the measure $\mu^*$ can
be explicitly identified. Throughout this section, we assume again
that $g(t) = 0$ for $t \leq 0$.

\medskip
\example{Example 2}  $N=1$ and $g$ is arbitrary.
\smallskip
This case is very easy since every measure is diffuse (recall that the only set of
zero capacity is the empty set). Hence, by Corollary 2, every measure
is good.  Thus, $\mu^* = \mu$ for every $\mu$.
\endexample
\medskip

\example{Example 3}  $N\geq 2$ and $g(t) = t^p$, $t\geq 0$, with
$1< p < \frac{N}{N-2}$.
\medskip

In this case, we have
again $\mu^* = \mu$ since, for {\it every}\/ measure $\mu$,
problem (0.1) admits a solution. This result was originally
established in 1975 by Ph.~B\'enilan and H.~Brezis (see [BB,
Appendix A], [B1], [B2], [B3], [B4] and also Theorem~B.5 below). The
crucial ingredient is
the compactness of the
imbedding of the space $\big\{u\in W^{1,1}_0 \, ; \, \Delta u \in \Cal
M\big\}$, equipped
with the norm $\|u\|_{W^{1,1}} + \|\Delta u\|_{\Cal M}$, into $L^q$ for
every $q < \frac{N}{N-2}$ (see Theorem~B.1 below).
\endexample
\medskip

\example{Example 4}  $N\geq 3$ and $g(t) = t^p$, $t \geq 0$, with $p\geq
\frac{N}{N-2}$.
\medskip

In this case, we have

\proclaim{Theorem 16}  For every measure $\mu$, we have
$$
\mu^* = \mu - (\mu_2)^+, \tag 8.1
$$
where $\mu = \mu_1+ \mu_2$ is the unique decomposition  of $\mu$ (in the sense
of Lemma~A.1) relative to the $W^{2, p'}$-capacity.
\endproclaim

\demo{Proof}  By a result of Baras-Pierre~[BP] (already mentioned in
the Introduction) we know that a measure $\nu \geq 0$ is a good
measure if and only if $\nu$ is diffuse with respect to the $W^{2,
p'}$-capacity.

Set
$$
\tilde \mu = \mu - (\mu_2)^+ = \mu_1- (\mu_2)^- \quad \text{and} \quad
\tilde \nu = (\mu_2)^+.\tag 8.2
$$
We claim that
$$
(\tilde \mu)^* = \tilde \mu  \quad \text{and} \quad (\tilde \nu)^* =
0. \tag 8.3
$$
Clearly, $(\tilde \mu)^+ = (\mu_1)^+$.
From the result of Baras-Pierre~[BP], we infer that $(\mu_1)^+$ is a
good measure. By Theorem~4, $\tilde\mu$ is also a good measure. Thus,
$(\tilde \mu)^* = \tilde \mu$. Since $\tilde\nu$ is a nonnegative
measure concentrated on a set of zero $W^{2,p}$-capacity, it follows
from [BP] that $(\tilde\nu)^* \leq 0$. Since $(\tilde\nu)^* \geq 0$,
we conclude that (8.3) holds.

\nd
Applying Theorem~8, we get
$$
\mu^* = (\tilde \mu + \tilde \nu)^* = (\tilde \mu)^* + (\tilde \nu)^*
= \tilde \mu = \mu - (\mu_2)^+,
$$
which is precisely (8.1).
\enddemo

\remark{Remark 9} In this example we see that the measure $\mu -
\mu^*$ is concentrated on a set $\Sigma$ whose $W^{2, p'}$-capacity
is zero. This is a better information than the general fact that $\mu
- \mu^*$ is concentrated on a set $\Sigma$ whose $H^1$-capacity is
zero.
\endremark
\endexample

\example{Example 5}  $N=2$ and $g(t) = \e^t-1$, $t \geq 0$.
\medskip
In this case, the identification of $\mu^*$ relies heavily on a result
of V\'azquez~[Va].

\proclaim{Theorem 17}  Given any measure $\mu$, let
$$
\mu = \mu_1 + \mu_2
$$
where $\mu_2$ is the purely atomic part of $\mu$ (this
corresponds to the decomposition of $\mu$ in the sense of Lemma~A.1,
where $\Cal Z$ consists of
countable sets). Write
$$
\mu_2 = \sum_i \alpha_i \delta_{a_i} \tag 8.4
$$
with $a_i \in \Omega$ distinct,  and $\sum|\alpha_i|<\infty$.
Then
$$
\mu^* = \mu - \sum_i (\alpha_i - 4\pi)^+ \delta_{a_i}.\tag 8.5
$$
\endproclaim

\demo{Proof}  By a result of V\'azquez~[Va], we know that a measure
$\nu$ is a good measure if and only if $\nu(\{x\})\leq 4\pi$ for every
$x\in\Omega$. (The paper of V\'azquez deals with the equation (0.1) in
all of $\Bbb R^2$ but the conclusion, and the proof, are the same for
a bounded domain).

Clearly, $\mu_1 (\{x\}) = 0$, $\forall x \in \Omega$.
From the result of V\'azquez~[Va] we infer that $\mu_1$ is a good
measure. Thus,
$$
(\mu_1)^* = \mu_1. \tag 8.6
$$
Let $a \in \Omega$ and $\alpha \in \Bbb R$. It is easy to see from [Va] that
$$
(\alpha \delta_a)^* = \min{\{\alpha, 4 \pi \}} \, \delta_a. \tag 8.7
$$
An induction argument applied to Theorem~8 and the
continuity of the mapping $\mu \mapsto \mu^*$ show that
$$
\mu^* = (\mu_1)^* + (\mu_2)^* = (\mu_1)^* + \sum_i{(\alpha_i
\delta_{a_i})^* }. \tag 8.8
$$
By (8.6)--(8.8), we have
$$
\mu^* = \mu_1 +  \sum_i{\min{\{\alpha_i, 4 \pi \}} \,
\delta_{a_i}} = \mu - \sum_i (\alpha_i - 4\pi)^+ \delta_{a_i}. \tag 8.9
$$
This establishes (8.5).
\enddemo

\endexample

We conclude this section with two interesting questions:

\medskip
\definition{Open problem 3}  Let $N=2$ and $g(t) = (\e^{t^2}-1)$, $t\geq
0$.  Is there an explicit formula for $\mu^*$~?
\enddefinition

\definition{Open problem 4}  Let $N\geq 3$ and $g(t) = (\e^t -1)$, $t
\geq 0$.  Is there an explicit formula for $\mu^*$~?
\enddefinition

A partial answer to Open problem~4 has been obtained by
Bartolucci-Leoni-Orsina-Ponce~[BLOP]. More precisely, they have
established the following:

\proclaim{Theorem~18}
Any measure $\mu$ such that $\mu \leq 4 \pi \Cal H^{N-2}$ is a good
measure.
\endproclaim

Here, $\Cal H^{N-2}$ denotes the $(N-2)$-Hausdorff measure. The
converse of Theorem~18 is not true. This was suggested by L.~V\'eron
in a personal communication; explicit examples are given in [P].
The characterization of good measures is still open; see however [MV6].

\bigskip


\subhead 9. Further directions and open problems \endsubhead
\medskip

\subhead 9.1. Vertical asymptotes \endsubhead
\medskip

Let $g: (-\infty,+1)\to \Bbb R$ be a continuous, nondecreasing
function such that $g(t)=0$, $\forall t \leq 0$, and such $ g(t)\to
+\infty$ as $ t \to +1$.
Let $(g_n)$ be a sequence of functions $g_n: \Bbb R\to \Bbb R$ which are
continuous, nondecreasing and satisfy the following conditions:
$$
\gather
0\leq g_1(t) \leq g_2(t)\leq \ldots\leq g(t)\quad\forall t<1, \tag 9.1\\
g_n(t) \to g(t)\quad \forall t<1 \quad \text{and}\quad g_n(t) \to
+\infty\quad \forall t \geq 1. \tag 9.2
\endgather
$$
If $N \geq 2$, then we also assume that
$$
\text{each $g_n$ has subcritical growth, i.e.,}\quad g_n(t)\leq C(|t|^p+
1)\quad \forall t \in \Bbb R, \tag 9.3
$$
for some constant $C$ and some $p<\frac{N}{N-2}$, possibly depending on $n$.

Given $\mu \in \Cal M(\Omega)$, let $u_n$ be a solution of
(0.10). Then $u_n\downarrow u^*$ in $\Omega$ as $n\uparrow
\infty$. Moreover (0.11) and (0.12) hold. We may therefore define
$\mu^* \in \Cal M(\Omega)$ by (0.13).

\definition{Open problem 5}
Study the properties of $u^*$ and the reduced measure $\mu^*$.
\enddefinition

Clearly, $u^*$ is the largest subsolution. But there are some major
differences in this case. When $N \geq 2$, Dupaigne-Ponce-Porretta~[DPP] have shown
that for any such $g$ one can find a nonnegative measure $\mu$ for
which the set $\{ \nu \in \Cal G \; ; \; \nu \leq \mu \}$ has {\it
no}\/ largest element. In particular, for such measure $\mu$, the
reduced measure $\mu^*$ cannot be the largest good measure $\leq
\mu$. They have also proved that the set of good measures $\Cal G$ is {\it not}\/ convex
for any $g$. We refer the reader to [DPP] for other results.

\medskip
Similar questions arise when $g$ is a multivalued graph. For example,
$$
g(r) =
\left\{
\alignedat0
0          & \quad \text{if } r < 1, \\
[0,\infty) & \quad \text{if } r = 1, \\
\emptyset  & \quad \text{if } r>1.
\endalignedat
\right.
$$
This is a simple model of one-sided  variational
inequality. The objective is to solve in some natural ``weak" sense
the multivalued equation
$$
\left\{
\alignedat2
-\Delta u + g(u) & \ni \mu && \quad \text{in }  \Omega,\\
u & = 0 && \quad \text{on } \partial \Omega,
\endalignedat
\right.
$$
for any given $\mu \in \Cal M(\Omega)$. This problem has been recently
studied by Brezis-Ponce~[BP4]. There were some partial results;
see, e.g., Baxter [Ba], Dall'Aglio-Dal Maso [DD],
Orsina-Prignet [OP], Brezis-Serfaty [BSe], and the references therein.

\bigskip


\subhead 9.2. Nonlinearities involving $\bs \nabla \bs u$ \endsubhead
\medskip

Consider the model problem:
$$
\left\{
\alignedat2
-\Delta u + u|\nabla u|^2 &= \mu && \quad \text{in }  \Omega,\\
u & = 0 && \quad \text{on } \partial \Omega,
\endalignedat
\right. \tag 9.4
$$
where $\mu \in \Cal M(\Omega)$. Problems of this type have been
extensively studied and it is known that they bear some similarities
with the problems discussed in this paper. In particular, it has been
proved in [BGO2] that (9.4) admits a solution
if and only if the measure $\mu$ is diffuse, i.e., $|\mu|(A) = 0$ for
every Borel set $A\subset\Omega$  such that $\capt{(A)}=0$. Moreover,
the solution is unique (see [BM]). When
$\mu$ is a general measure, not necessarily diffuse, it would be
interesting to apply to (9.4) the same strategy as in this paper. More
precisely, to prove that approximate solutions converge to the solution
of (9.4), where $\mu$ is replaced by its diffuse part $ \mu_{\trd}$
(in the sense of Lemma A.1, relative to the Borel sets whose
$H^1$-capacity are zero):
$$
\left\{
\alignedat2
-\Delta u + u|\nabla u|^2 & = \mu_{\trd} && \quad \text{in }  \Omega,\\
u & = 0 && \quad \text{on } \partial \Omega,
\endalignedat
\right. \tag 9.5
$$
which possesses a unique solution. There are several ``natural''
approximations. For example, one may truncate the nonlinearity
$g(u,\nabla u)=u|\nabla u|^2$ and replace it by
$g_n(u,\nabla u)= \frac{n}{n+ |g(u,\nabla u)|} \, g(u,\nabla u)$.
It is easy to see (via a
Schauder fixed point argument in $W^{1,1}_0$) that the corresponding
equation
$$
\left\{
\alignedat2
-\Delta u_n + g_n(u_n,\nabla u_n) & = \mu && \quad \text{in }  \Omega,\\
u_n & = 0 && \quad \text{on } \partial \Omega,
\endalignedat
\right. \tag 9.6
$$
admits a solution $u_n$.

\definition{Open problem 6}
Is it true that $(u_n)$ converges to the
solution of (9.5)?
\enddefinition

Another possible approximation consists of smoothing $\mu$: let $u_n$
be a solution of
$$
\left\{
\alignedat2
-\Delta u_n + u_n |\nabla u_n|^2 & = \mu_n && \quad \text{in }  \Omega,\\
u_n & = 0 && \quad \text{on } \partial \Omega,
\endalignedat
\right. \tag 9.7
$$
where $\mu_n = \rho_n * \mu$, as in Section~4. It has been proved by
Porretta~[Po] that if $\mu \geq 0$, then $u_n \to u$ in $L^1(\Omega)$,
where $u$ is the solution of (9.5). We have been informed by
A.~Porretta that the same conclusion holds for any measure $\mu$, by
using a substantial modification of the argument in [Po].

\bigskip


\subhead 9.3. Measures as boundary data \endsubhead
\medskip

Consider the problem
$$
\cases
-\Delta u + g(u) = 0 & \text{in }  \Omega,\\
\qquad\qquad\;\;\,  u = \mu &\text{on } \partial \Omega,
\endcases\tag 9.8
$$
where $\mu$ is a measure on $\partial \Omega$ and $g: \Bbb R \to
\Bbb R$ is a continuous, nondecreasing function satisfying (0.5). It
has been proved by H.$ $ Brezis (1972, unpublished) that (9.8) admits a
unique weak solution when $\mu$ is any $L^1$ function (for a general
nonlinearity $g$). When $g$ is a power, the study of (9.8) for
measures was initiated by
Gmira-V\'eron~[GV], and has vastly expanded in recent years; see the papers
of Marcus-V\'eron~[MV1], [MV2], [MV3], [MV4]. Important motivations coming from the
theory of probability --- and the use of probabilistic methods --- have
reinvigorated the whole subject; see the pioneering papers of
Le Gall~[LG1], [LG2], the recent books of Dynkin~[D1], [D2], and the numerous
references therein. It is known that (9.8) has no solution if
$g(t)=t^p$, $t\geq 0$, with $p \geq p_\trc= \frac{N+1}{N-1}$ and $\mu = \delta_a$, $a
\in \partial \Omega$ (see [GV]). Therefore, it is interesting to develop for
(9.8) the same program as in this paper.
More precisely, let $(g_k)$ be a sequence of functions $g_k : \Bbb R \to \Bbb R$ which
are continuous, nondecreasing, and satisfy (0.7) and (0.8). Assume
in addition that each $g_k$ is, e.g., bounded. Then, for every $\mu
\in \Cal M(\partial\Omega)$, there exists a unique solution $u_k$ of
$$
\left\{
\alignedat2
-\Delta u_k + g_k(u_k) &= 0 && \quad \text{in }  \Omega,\\
u_k &= \mu && \quad \text{on } \partial \Omega,
\endalignedat
\right. \tag 9.9
$$
in the sense that $u_k \in L^1(\Omega)$ and
$$
- \int_\Omega u_k \Delta \zeta + \int_\Omega g_k(u_k) \zeta =
-  \int_{\partial\Omega}  \frac{\partial \zeta}{\partial n} \, d\mu \quad
\forall \zeta \in C_0^2(\overline\Omega), \tag 9.10
$$
where $\frac{\partial}{\partial n}$ denotes the derivative with
respect to the outward normal of $\partial\Omega$.

\noindent
We have the following:

\proclaim{Theorem 19}
As $k \uparrow \infty$, $u_k \downarrow u^*$ in $L^1(\Omega)$,
where $u^*$ satisfies
$$
\left\{ \alignedat2
-\Delta u^* + g(u^*) & = 0 && \quad \text{in }  \Omega,\\
u^* & = \mu^* && \quad \text{on } \partial \Omega,
\endalignedat \right.\tag 9.11
$$
for some $\mu^* \in \Cal M(\partial\Omega)$ such that $\mu^* \leq
\mu$. More precisely, $g(u^*) \rho_0 \in L^1(\Omega)$, where $\rho_0(x) = d(x,
\partial\Omega)$, and
$$
- \int_\Omega u^* \Delta \zeta + \int_\Omega g(u^*) \zeta =
-  \int_{\partial\Omega}  \frac{\partial \zeta}{\partial n} \, d\mu^* \quad
\forall \zeta \in C_0^2(\overline\Omega). \tag 9.12
$$
In addition, $u^*$ is the largest subsolution of\/ \rom{(9.8)}, i.e.,
if $v \in L^1(\Omega)$ is any function satisfying $g(v) \rho_0 \in L^1(\Omega)$ and
$$
- \int_\Omega v \Delta \zeta + \int_\Omega g(v) \zeta \leq
-  \int_{\partial\Omega}  \frac{\partial \zeta}{\partial n} \, d\mu \quad
\forall \zeta \in C_0^2(\overline\Omega), \; \zeta \geq 0 \text{ in }
\Omega, \tag 9.13
$$
then $v \leq u^*$ a.e.$ $ in $\Omega$.
\endproclaim

\demo{Proof}
By comparison (see Corollary~B.2), we know that $(u_k)$ is
non-increasing. By standard estimates, we have
$$
\int_\Omega |u_k| + \int_\Omega g_k(u_k) \rho_0 \leq C
\|\mu\|_{\Cal M(\partial\Omega)} \quad \forall k \geq 1.
$$
In addition, (see [B5, Theorem~3])
$$
\|u_k\|_{C^1(\overline\omega)} \leq C_\omega \quad \forall k \geq 1,
$$
for every $\omega \subset\subset \Omega$.
Thus, $u_k$ converges in $L^1(\Omega)$ to a limit, say $u^*$.
Moreover,
$$
g_k(u_k) \to g(u^*) \quad \text{in } L^\infty\loc(\Omega).
$$
Let $\zeta_0 \in C_0^2(\overline\Omega)$ be the solution of
$$
\left\{ \alignedat2
- \Delta \zeta_0 & = 1 && \quad \text{in } \Omega,\\
\zeta_0 & = 0          && \quad \text{on } \partial\Omega.
\endalignedat \right.
$$
Since $\big(g_k(u_k) \zeta_0\big)$ is uniformly bounded in $L^1(\Omega)$, then
up to a subsequence
$$
g_k(u_k)\zeta_0 \overset{*}\to\rightharpoonup g(u^*)\zeta_0 + \lambda \quad
\text{in } \big[C(\overline\Omega)\big]^*, \tag 9.14
$$
for some $\lambda \in \Cal M(\partial\Omega)$, $\lambda \geq 0$. We claim that
$$
\int_\Omega g_k(u_k)\zeta \to \int_\Omega g(u^*)\zeta +
\int_{\partial\Omega} \frac{\partial \zeta}{\partial n}
\frac{1}{\frac{\partial\zeta_0}{\partial n}} \, d\lambda \quad
\forall \zeta \in C_0^2(\overline\Omega). \tag 9.15
$$
In fact, given $\zeta \in C_0^2(\overline\Omega)$, define $\gamma =
\zeta/\zeta_0$. It is easy to see that $\gamma \in
C(\overline\Omega)$ and $\gamma = \frac{\partial \zeta}{\partial n}
\frac{1}{\frac{\partial\zeta_0}{\partial n}}$ on
$\partial\Omega$. Using $\gamma$ as a test function in (9.14), we obtain
(9.15).

\noindent
Let $k \to \infty$ in (9.10). In view of (9.15), we conclude that $u^*$ satisfies
(9.12), where $\mu^*$ is given by
$$
\mu^* = \mu + \frac{1}{\frac{\partial\zeta_0}{\partial n}}
\,\lambda \leq \mu.
$$
Finally, it follows from Corollary B.2 that if $v$ is a subsolution of
(9.8), then $v\leq u_k$ a.e., $\forall k \geq 1$, and thus $v\leq u^*$
a.e.
\enddemo

Some natural questions have been addressed and the following results
will be presented in a forthcoming paper (see [BP3]):
\medskip

{\rm(a)} the reduced measure $\mu^*$ is the largest good measure $\leq \mu$; in other
words, if $\nu \in \Cal M(\partial\Omega)$ is a good measure (i.e., (9.8) has
a solution with boundary data $\nu$) and if $\nu \leq \mu$, then $\nu
\leq \mu^*$;
\medskip

{\rm(b)}  $\mu-\mu^*$ is concentrated on a subset of $\partial\Omega$ of zero $\Cal
H^{N-1}$-measure (i.e., $(N-1)$-dimensional Lebesgue measure on
$\partial\Omega$) and this fact is ``optimal'', in the sense that any
measure $\nu \geq 0$ which is singular with respect to $\Cal
H^{N-1}\lfloor_{\partial\Omega}$ can be written as $\nu = \mu - \mu^*$
for some $\mu \geq 0$ and some $g$;
\medskip

{\rm(c)} if $\mu$ is a measure on $\partial\Omega$ which is good for
every $g$, then $\mu^+ \in L^1(\partial\Omega)$;
\medskip

{\rm(d)} given any $g$, there exists some measure $\mu \geq 0$ on
$\partial\Omega$ which is good for $g$, but $\mu \not\in L^1(\partial\Omega)$.
\medskip

\noindent
When $g(t) = t^p$, $t \geq 0$, with $p \geq \frac{N+1}{N-1}$, a known
result (see, e.g., [MV3]) asserts that $\mu \in \Cal M(\partial\Omega)$ is a
good measure if and only if $\mu^+(A) = 0$ for every Borel set $A
\subset \partial\Omega$ such that $C_{2/p, p'}(A) = 0$, where $C_{2/p,
p'}$ refers to the Bessel capacity on $\partial\Omega$. In this case, we have
\medskip

{\rm(e)} the reduced measure $\mu^*$ is given by $\mu^* = \mu -
(\mu_2)^+$, where $\mu = \mu_1 + \mu_2$ is the
decomposition of $\mu$, in the sense of Lemma~A.1, relative to $C_{2/p, p'}$.
\medskip

\noindent
In contrast with Example~5, we do not know what the reduced measure
$\mu^*$ is when $N=2$ and $g(t) = \e^t -1$, $t \geq 0$.

\bigskip
Similar issues can be investigated for the parabolic equations
$$
\left\{ \alignedat1
u_t - \Delta u + g(u) & = \mu,\\
u(0) & =0,
\endalignedat \right.
\qquad \text{or} \qquad
\left\{ \alignedat1
u_t - \Delta u + g(u) & = 0,\\
u(0) & =\mu.
\endalignedat \right.
$$
\smallskip

\bigskip


\subhead Appendix A: Decomposition of measures into diffuse and concentrated parts
\endsubhead

\medskip
The following result is taken from [FST]. We reproduce their proof for the convenience of the reader.
\medskip

\proclaim{Lemma A.1}  Let $\mu$ be a bounded Borel measure in $\Bbb
R^N$ and let $\Cal Z$ be a collection of Borel sets such that:
\medskip

{\rm (a)} $\Cal Z$ is closed with respect to finite or countable unions;

\medskip
{\rm (b)} $A\in \Cal Z$ and $A'\subset A$ Borel\/ $\Rightarrow$ $A' \in \Cal Z$.

\medskip
\noindent Then $\mu$ can be represented in the form
$$
\mu= \mu_1 + \mu_2, \tag A.1
$$
where $\mu_1$ and $\mu_2$ are bounded Borel measures such that
$$
\mu_1(A) = 0 \quad \forall A\in \Cal Z \quad \text{and}  \quad  \mu_2 \text{ vanishes outside a set
} A_0 \in \Cal Z.
$$
This representation is unique.
\endproclaim

\demo{Proof} First assume that $\mu$ is nonnegative.  Denote
$$
X_\mu = \sup{\big\{ \mu(A) \, ; \, A \in \Cal Z \big\}}. 
$$
Let $\{A_n\}$ be an increasing sequence of sets in $\Cal Z$ such that
$$
\mu(A_n) \to X_\mu. 
$$
Let $A_0= \bigcup_n A_n$ and put
$$
\mu_1(B) = \mu(B\cap A^{\trc}_0), \quad \mu_2(B) = \mu(B\cap A_0),
$$
for every Borel set $B$.  Since $A_0\in \Cal Z$, it remains to verify that
$\mu_1$ vanishes on sets of $\Cal Z$.  By contradiction, suppose that there
exists $E \in \Cal Z$ such that $\mu_1(E)> 0$.  Let $E_1=E \cap A^{\trc}_0$.
Then $\mu(E_1) > 0$ and $E_1 \in \Cal Z$.  It follows that $A_0\cup E_1\in
\Cal Z$ and $\mu(A_0 \cup E_1) > X_\mu$.  Contradiction.

\nd
If $\mu$ is a signed measure, apply the above to $\mu^+$ and $\mu^-$.
The uniqueness is obvious.
\enddemo

\bigskip


\subhead Appendix B: Standard existence, uniqueness and comparison results \endsubhead
\medskip

In this appendix, we collect some well-known results (and a few new
ones) which are used
throughout this paper. For the convenience of the reader, we shall
sketch some of the proofs.

We start with the existence, uniqueness and regularity of solutions of
the linear problem
$$
\left\{ \alignedat2
-\Delta u &= \mu  && \quad \text{in } \Omega,\\
u & =0  && \quad \text{on } \partial\Omega,
\endalignedat \right. \tag B.1
$$
where $\mu \in \Cal M(\Omega)$.

\proclaim{Theorem B.1}
Given $\mu \in \Cal M(\Omega)$, there exists a unique $u \in
L^1(\Omega)$ satisfying
$$
- \int_\Omega u \Delta \zeta = \int_\Omega \zeta \, d\mu \quad \forall
  \zeta \in C_0^2(\overline\Omega). \tag B.2
$$
Moreover, $u \in W_0^{1,q}(\Omega)$ for every $1 \leq q <
\frac{N}{N-1}$, with the estimates
$$
\|u\|_{L^{q^*}} \leq C \|\nabla u\|_{L^q} \leq C \|\mu\|_{\Cal M},
\tag B.3
$$
where $\frac{1}{q^*} = \frac{1}{q} - \frac{1}{N}$. In particular,  $u
\in L^p(\Omega)$ for every $1 \leq p < \frac{N}{N-2}$, and $u$
satisfies
$$
\int_\Omega \nabla u \cdot \nabla \psi = \int_\Omega \psi \, d\mu
\quad \forall \psi \in W_0^{1,r}(\Omega), \tag B.4
$$
for any $r > N$.
\endproclaim

The proof of Theorem~B.1 relies on a standard duality argument and
shall be omitted; see [S, Th\'eor\`eme~8.1].

\smallskip
We now establish a weak form of the maximum principle:

\proclaim{Proposition B.1}
Let $v \in W_0^{1,1}(\Omega)$ be such that
$$
- \int_\Omega v \Delta \varphi \leq 0 \quad \forall \varphi \in
  C_\trc^\infty(\Omega), \; \varphi \geq 0 \text{ in } \Omega. \tag B.5
$$
Then
$$
- \int_\Omega v \Delta \zeta \leq 0 \quad \forall \zeta \in
  C_0^2(\overline\Omega), \; \zeta \geq 0 \text{ in } \Omega \tag B.6
$$
and, consequently,
$$
v \leq 0 \quad \text{a.e.} \tag B.7
$$
\endproclaim

\demo{Proof}
From (B.5) we have
$$
\int_\Omega \nabla v \cdot \nabla \varphi \leq 0 \quad \forall \varphi \in
  C_\trc^\infty(\Omega), \; \varphi \geq 0 \text{ in } \Omega
$$
so that, by density of $C_\trc^\infty(\Omega)$ in $C_\trc^2(\Omega)$,
$$
\int_\Omega \nabla v \cdot \nabla\varphi \leq 0 \quad \forall \varphi \in
  C_\trc^2(\Omega), \; \varphi \geq 0 \text{ in } \Omega.
$$
Let $(\gamma_n)$ be a sequence in $C_\trc^\infty(\Omega)$ such that $0
\leq \gamma_n \leq 1$, $\gamma_n(x)  = 1$ if $d(x, \partial\Omega) >
\frac{1}{n}$, and $|\nabla\zeta_n| \leq Cn$, $\forall n \geq 1$. For
any $\zeta \in C_0^2(\overline\Omega)$, $\zeta \geq 0$, we have
$$
\int_\Omega \nabla v \cdot (\gamma_n \nabla \zeta +  \zeta
\nabla \gamma_n) = \int_\Omega \nabla v \cdot
\nabla(\gamma_n \zeta) \leq 0
. \tag B.8
$$
Note that
$$
\int_\Omega |\nabla v| |\nabla \gamma_n| \zeta \leq C n \underset{d(x,
\partial\Omega) \leq \frac{1}{n}}\to\int |\nabla v| \zeta \leq C \underset{d(x,
\partial\Omega) \leq \frac{1}{n}}\to\int |\nabla v| \to 0 \quad
\text{as } n \to \infty.
$$
Thus, as $n \to \infty$ in (B.8), we obtain
$$
\int_\Omega \nabla v \cdot \nabla\zeta \leq 0 \quad \forall \zeta \in
  C_0^2(\overline\Omega), \; \zeta \geq 0 \text{ in } \Omega,
$$
which yields (B.6) since $v \in W_0^{1,1}(\Omega)$. Inequality (B.7)
is a trivial consequence of (B.6).
\enddemo

\proclaim{Lemma B.1}
Let $p : \Bbb R \to \Bbb R$, $p(0) = 0$, be a bounded nondecreasing
continuous function. Given $f \in L^1(\Omega)$, let $u \in
L^1(\Omega)$ be the unique solution of
$$
- \int_\Omega u \Delta \zeta = \int_\Omega f \zeta \quad \forall \zeta
  \in C_0^2(\overline\Omega). \tag B.9
$$
Then
$$
\int_\Omega f p(u) \geq 0. \tag B.10
$$
\endproclaim

\demo{Proof}
Clearly, it suffices to establish the lemma for $p \in C^2(\Bbb R)$.
Assume for the moment $f \in C^\infty(\overline\Omega)$. In this case,
$u \in C_0^2(\overline\Omega)$. Since $p(0) = 0$, we have $p(u) \in
C_0^2(\overline\Omega)$. Using $p(u)$ as a test function in (B.9), we
get
$$
\int_\Omega f p(u) = \int_\Omega p'(u) |\nabla u|^2 \geq 0.
$$
This establishes the lemma for $f$ smooth. The general case when $f$
is just an $L^1$-function, not necessarily smooth, easily follows by
density.
\enddemo

\proclaim{Proposition B.2}
Given $f \in L^1(\Omega)$, let $u$ be the unique solution of
\rom{(B.9)}. Then, for every $M > 0$, we have
$$
\underset{[u \geq M]}\to\int f \geq 0 \quad \text{and} \quad
\underset{[u \leq -M]}\to\int f \leq 0 . \tag B.11
$$
In particular,
$$
\underset{[|u| \geq M]}\to\int f \sign{(u)} \geq 0. \tag B.12
$$
\endproclaim

Above, we denote by $\sign$ the function $\sign{(t)} = 1$ if $t >0$, $\sign{(t)}
= -1$ if $t < 0$, and $\sign{(0)} = 0$.

\demo{Proof}
Clearly, it suffices to establish the first inequality in (B.11). Let
$(p_n)$ be a sequence of continuous functions in $\Bbb R$ such that each
$p_n$ is nondecreasing, $p_n(t) = 1$ if $t \geq M$ and $p_n(t) = 0$ if $t
\leq M - \frac{1}{n}$. By the previous lemma,
$$
\int_\Omega f p_n(u) \geq 0 \quad \forall n \geq 1.
$$
As $n \to \infty$, the result follows.
\enddemo

\proclaim{Proposition B.3}
Let $v \in L^1(\Omega)$, $f \in L^1(\Omega)$ and $\nu \in \Cal
M(\Omega)$ satisfy
$$
-\int_\Omega v \Delta \zeta + \int_\Omega f \zeta = \int_\Omega \zeta
 \, d\nu \quad \forall \zeta \in C_0^2(\overline\Omega).\tag B.13
$$
Then
$$
\underset{[v > 0]}\to\int f \leq \|\nu^+\|_{\Cal M} \tag B.14
$$
and thus
$$
\int_\Omega f \sign{(v)} \leq \|\nu\|_{\Cal M}. \tag B.15
$$
\endproclaim

\demo{Proof}
Let $\nu_n = \rho_n * \nu$ (here we use the same notation as in
Section~4). Let $v_n$ denote the solution of (B.13) with $\nu$
replaced by $\nu_n$. By Lemma~B.1, we have
$$
\int_\Omega (\nu_n - f) \,p(v_n) \geq 0,
$$
where $p$ is any function satisfying the assumptions of the lemma.
Thus, if $0 \leq p(t) \leq 1$, $\forall t \in \Bbb R$, then we have
$$
\int_\Omega f  p(v_n) \leq \int_\Omega \nu_n p(v_n) \leq \int_\Omega
(\nu_n)^+ \leq \|\nu^+\|_{\Cal M}.
$$
Let $n \to \infty$ to get
$$
\int_\Omega f  p(v)  \leq \|\nu^+\|_{\Cal M}. \tag B.16
$$
Apply (B.16) to a sequence of nondecreasing continuous functions
$(p_n)$ such that $p_n(t) = 0$ if $t \leq 0$ and $p_n(t) = 1$ if $t
\geq \frac{1}{n}$. As $n \to \infty$, we obtain (B.14).
\enddemo

\smallskip
An easy consequence of Proposition~B.3 is the following

\proclaim{Corollary B.1}
Let $g: \Bbb R \to \Bbb R$ be a continuous, nondecreasing function
such that $g(0)= 0$. Given $\mu \in \Cal M(\Omega)$, then the equation
$$
\left\{ \alignedat2
-\Delta u + g(u) & = \mu && \quad \text{in }  \Omega,\\
u & = 0 && \quad \text{on } \partial \Omega,
\endalignedat \right.\tag B.17
$$
has\/ {\rm at most} one solution $u \in L^1(\Omega)$ with $g(u) \in
L^1(\Omega)$. Moreover,
$$
\int_\Omega |g(u)| \leq \|\mu\|_{\Cal M} \quad \text{and}
\quad \int_\Omega |\Delta u| \leq 2 \|\mu\|_{\Cal M}. \tag B.18
$$
If\/ \rom{(B.17)} has a solution for $\mu_1, \mu_2 \in \Cal
M(\Omega)$, say $u_1, u_2$, resp., then
$$
\int_\Omega \big[ g(u_1) - g(u_2) \big]^+ \leq \big\| (\mu_1 -
\mu_2)^+ \big \|_{\Cal M}. \tag B.19
$$
In particular,
$$
\int_\Omega \big| g(u_1) - g(u_2) \big| \leq \| \mu_1 -
\mu_2 \|_{\Cal M}. \tag B.20
$$
\endproclaim

\smallskip

We now recall the following unpublished result of H.$ $ Brezis from
1972 (see, e.g., \cite{GV}):

\proclaim{Proposition B.4}
Given $f \in L^1(\Omega; \rho_0 \, dx)$ and $h \in L^1(\partial\Omega)$, there exists
a unique $u \in L^1(\Omega)$ such that
$$
- \int_\Omega u \Delta \zeta = \int_\Omega f \zeta - \int_{\partial\Omega} h \frac{\partial\zeta}{\partial n} \quad
\forall \zeta \in C_0^2(\overline\Omega). \tag B.21
$$
In addition, there exists $C > 0$ such that
$$
\|u\|_{L^1} \leq C \left( \|f \rho_0\|_{L^1(\Omega)} +
\|h\|_{L^1(\partial\Omega)} \right). \tag B.22
$$
\endproclaim

\smallskip

We now establish the following

\proclaim{Lemma B.2}
Given $f \in L^1(\Omega; \rho_0 \, dx)$, let $u \in L^1(\Omega)$ be the unique
solution of\/ \rom{(B.21)} with $h = 0$.
Then
$$
k \underset{d(x, \partial\Omega) < \frac{1}{k}}\to\int |u| \to 0 \quad
\text{as } k \to \infty. \tag B.23
$$
\endproclaim

\demo{Proof}

\noindent
{\it Step 1.}  Proof of the lemma when $f \geq 0$.
\smallskip

Since $f \geq 0$,  we have $u \geq 0$. Let
$H \in C^2(\Bbb R)$ be a nondecreasing concave function such that
$H(0) = 0$, $H''(t) =
- 1$ if $t \leq 1$ and $H(t) = 1$ if $t \geq 2$. We denote by $\zeta_0
\in C_0^2(\overline\Omega)$, $\zeta_0 \geq 0$, the solution of
$$
\left\{ \alignedat2
- \Delta \zeta_0 & = 1 && \quad \text{in } \Omega,\\
\zeta_0 & = 0          && \quad \text{on } \partial\Omega.
\endalignedat \right.
$$
For any $k \geq 1$, let $w_k = \frac{1}{k} H(k \zeta_0)$. By
construction, $w_k \in C_0^2(\overline\Omega)$ and
$$
\Delta w_k = k H''(k \zeta_0) |\nabla \zeta_0|^2 + H'(k \zeta_0)
\Delta \zeta_0 \leq - k \chi_{[\zeta_0 \leq \frac{1}{k}]} |\nabla
\zeta_0|^2.
$$
Thus,
$$
- \int_\Omega u \Delta w_k \geq k \underset{[\zeta_0 \leq
\frac{1}{k}]}\to\int |\nabla\zeta_0|^2 u . \tag B.24
$$
Use $w_k$ as a test function in (B.21) (recall that $h = 0$). It follows from (B.24) that
$$
k \underset{[\zeta_0 \leq \frac{1}{k}]}\to\int |\nabla\zeta_0|^2 u \leq \int_\Omega w_k f.
\tag B.25
$$
By Hopf's lemma, we have
$|\nabla \zeta_0|^2 \geq \alpha_0 > 0$ in some
neighborhood of $\partial\Omega$ in $\overline\Omega$. In particular,
there exists $c>0$ such that $c \zeta_0(x) \leq d(x, \partial\Omega) \leq \frac{1}{c}
\zeta_0(x)$ for all $x \in \overline\Omega$. Thus, for $k \geq 1$
sufficiently large, we have
$$
\alpha_0 k \underset{d(x, \partial\Omega) \leq \frac{c}{k}}\to\int |\nabla\zeta_0|^2 u \leq
\int_\Omega w_k f. \tag B.26
$$
Note that the right-hand side of (B.26) tends to 0 as $k \to \infty$.
In fact, we have $w_k \leq C \zeta_0$, $\forall k \geq 1$, and $w_k \leq \frac{1}{k} H(k \zeta_0) \to 0$
a.e. Thus, by dominated convergence,
$$
\int_\Omega w_k f \to 0 \quad \text{as } k \to \infty. \tag B.27
$$
Combining (B.26) and (B.27), we obtain (B.23).

\smallskip
\noindent
{\it Step 2.}  Proof of the lemma completed.
\smallskip

Let $v \in L^1(\Omega)$ denote the unique solution of
$$
- \int_\Omega v \Delta \zeta = \int_\Omega |f| \zeta  \quad \forall
  \zeta \in C_0^2(\overline\Omega). \tag B.28
$$
By comparison, we have $|u| \leq v$. On the other hand, $v$ satisfies
the assumption of Step~1. Thus,
$$
{k \underset{d(x, \partial\Omega) <
\frac{1}{k}}\to\int |u|} \leq k \underset{d(x,
\partial\Omega) < \frac{1}{k}}\to\int v \to 0 \quad \text{as } k \to \infty. \tag B.29
$$
This establishes Lemma B.2.
\enddemo

\medskip

The next result is a new variant of Kato's inequality, where the test
function $\zeta$ need not have compact support in $\Omega$:

\proclaim{Proposition B.5}
Let $u \in L^1(\Omega)$ and $f \in L^1(\Omega; \rho_0 \, dx)$ be such that
$$
- \int_\Omega u \Delta \zeta \leq \int_\Omega f \zeta \quad \forall
  \zeta  \in C_0^2(\overline\Omega), \; \zeta \geq 0 \text{ in } \Omega. \tag B.30
$$
Then
$$
- \int_\Omega u^+ \Delta \zeta \leq \underset{[u \geq 0]}\to\int f \zeta \quad \forall
  \zeta  \in C_0^2(\overline\Omega), \; \zeta \geq 0 \text{ in } \Omega. \tag B.31
$$
\endproclaim

\demo{Proof}
We first notice that
$$
- \int_\Omega u^+ \Delta \varphi \leq \underset{[u \geq 0]}\to\int f
 \varphi \quad \forall \varphi \in C_\trc^\infty(\Omega), \; \varphi \geq 0 \text{ in } \Omega. \tag B.32
$$
In fact, by (B.30) we have $- \Delta u \leq f$ in $\Cal D'(\Omega)$.
Then, Theorem 7 yields
$$
(- \Delta u^+)_\trd \leq \chi_{[u \geq 0]} (- \Delta u)_\trd \leq \chi_{[u \geq 0]} f \quad \text{and} \quad
(- \Delta u^+)_\trc = (- \Delta u)_\trc^+ \leq (f)_\trc^+ = 0.
$$
Thus,
$$
- \Delta u^+ = (- \Delta u^+)_\trd + (- \Delta u^+)_\trc \leq
  \chi_{[u \geq 0]} f \quad \text{in } \Cal D'(\Omega),
$$
which is precisely (B.32).

\smallskip
\nd
Let $(\gamma_k) \subset C_\trc^\infty(\Omega)$ be a sequence such
that $0 \leq \gamma_k \leq 1$ in $\Omega$, $\gamma_k(x) = 1$ if
$d(x, \partial\Omega) \geq \frac{1}{k}$, $\|\nabla\gamma_k \|_{L^\infty} \leq
k$, and $\|\Delta\gamma_k \|_{L^\infty} \leq C k^2$.
Given $\zeta \in C_0^2(\overline\Omega)$, $\zeta\geq 0$, we apply
(B.32) with $\varphi = \zeta \gamma_k$ to get
$$
- \int_\Omega u^+ \Delta ( \zeta \gamma_k) \leq \underset{[u \geq 0]}\to\int f
 \zeta \gamma_k.\tag B.33
$$
Consider again the unique solution $v \geq 0$ of (B.28). By comparison
we have $ u \leq v$ a.e. and thus $u^+ \leq v$ a.e.   From  Lemma~B.2 we see
that
$$
\int_\Omega u^+ |\nabla\zeta| |\nabla \gamma_k| \leq Ck
\underset{d(x, \partial\Omega) <
\frac{1}{k}}\to\int u^+ \to 0 \quad \text{as } k \to \infty. \tag B.34
$$
Similarly,
$$
\int_\Omega u^+ \zeta |\Delta \gamma_k| \leq Ck \underset{d(x, \partial\Omega) <
\frac{1}{k}}\to\int u^+ \to 0 \quad \text{as } k \to \infty. \tag B.35
$$
Let $k\to \infty$ in (B.33). Using (B.34) and (B.35), we obtain (B.31).
\enddemo

\medskip

\remark{Remark B.1}
There is an alternative proof of Proposition~B.5. First, one shows
that (B.30) implies that there exist two measures $\mu \leq 0$, $\lambda
\leq 0$, where $\mu \in \Cal M(\partial\Omega)$ and $\lambda$ is locally
bounded in $\Omega$, with $\int_\Omega \rho_0 \, d|\lambda| < \infty$,
satisfying
$$
- \int_\Omega u \Delta\zeta = \int_\Omega f \zeta + \int_\Omega \zeta \,
  d\lambda - \int_{\partial\Omega} \frac{\partial\zeta}{\partial n} \,
  d\mu \quad \forall \zeta \in C_0^2(\overline\Omega). \tag B.36
$$
[The existence of $\lambda$ is fairly straightforward, and the existence
of $\mu$ is a consequence of Herglotz's theorem concerning positive
superharmonic functions].

\nd
Then, inequality (B.31) follows from (B.36) using the same strategy as in the
proof of Lemma~1.5 in [MV2].
\endremark

\medskip

As a consequence of Proposition~B.5, we have the following

\proclaim{Corollary B.2}
Let $g_1, g_2 : \Bbb R \to \Bbb R$ be two continuous nondecreasing
functions such that $g_1 \leq g_2$.
Let $u_k \in L^1(\Omega)$, $k = 1, 2$, be such that $g_k(u_k) \in
L^1(\Omega; \rho_0 \, dx)$. If
$$
- \int_\Omega (u_2 - u_1) \Delta \zeta + \int_\Omega \big[g_2(u_2) -
g_1(u_1)\big] \zeta \leq 0 \quad \forall \zeta \in C_0^2(\overline\Omega), \;
\zeta \geq 0 \text{ in } \Omega, \tag B.37
$$
then
$$
u_2 \leq u_1 \quad \text{a.e.} \tag B.38
$$
\endproclaim

\demo{Proof}
Applying Proposition~B.5 to $u = u_2 - u_1$ and $f = g_1(u_1) -
g_2(u_2)$, we have
$$
- \int_\Omega (u_2 - u_1)^+ \Delta \zeta \leq - \int_\Omega \big[ g_2(u_2) - g_1(u_1)
\big]^+ \zeta \leq 0 \quad \forall
  \zeta  \in C_0^2(\overline\Omega), \; \zeta \geq 0 \text{ in } \Omega.
$$
This immediately implies that $u_2 \leq u_1$ a.e.
\enddemo

\smallskip
We now present some general existence results for problem
(B.17). Below, $g: \Bbb R \to \Bbb R$ denotes a continuous, nondecreasing function,
such that $g(0)= 0$.

\proclaim{Theorem B.2 (Brezis-Strauss~[BS])}
For every $f \in L^1(\Omega)$, the equation
$$
\left\{ \alignedat2
-\Delta u + g(u) & = f && \quad \text{in }  \Omega,\\
u & = 0 && \quad \text{on } \partial \Omega,
\endalignedat \right.\tag B.39
$$
has a unique solution $u \in L^1(\Omega)$ with $g(u) \in L^1(\Omega)$.
\endproclaim

\demo{Proof}
We first observe that if $f \in C^\infty(\overline\Omega)$, then
(B.39) always has a solution $u \in C^1(\overline\Omega)$ (easily
obtained via minimization).

\nd
For a general $f \in L^1(\Omega)$, let $(f_n)$ be a sequence of smooth
functions on $\overline\Omega$, converging to $f$ in $L^1(\Omega)$. For
each $f_n$, let $u_n$ denote the corresponding solution of (B.39). By
(B.20), the sequence $(g(u_n))$ is Cauchy in $L^1(\Omega)$. We then
conclude from (B.3) that $(u_n)$ is also Cauchy in
$L^1(\Omega)$, so that
$$
u_n \to u \quad \text{and} \quad g(u_n) \to g(u) \quad \text{in } L^1(\Omega).
$$
Thus $u$ is a solution of (B.39). The uniqueness follows
from Corollary~B.1.
\enddemo

\proclaim{Theorem B.3 (Brezis-Browder~[BBr])}
For every $T \in H^{-1}(\Omega)$, the equation
$$
\left\{ \alignedat2
-\Delta u + g(u) & = T && \quad \text{in }  \Omega,\\
u & = 0 && \quad \text{on } \partial \Omega,
\endalignedat \right.\tag B.40
$$
has a unique solution $u \in H_0^1(\Omega)$ with $g(u) \in L^1(\Omega)$.
\endproclaim

\demo{Proof}
Assume $g$ is uniformly bounded. In this case, the existence of $u$
presents no difficulty, e.g., via a minimization argument in
$H_0^1(\Omega)$. In particular, we see that $u \in H_0^1(\Omega)$.

\nd
For a general nonlinearity $g$, let $(g_n)$ be the sequence given by
$g_n(t) = g(t)$ if $|t| \leq n$, $g_n(t) = g(n)$ if $t > n$, and
$g_n(t) = g(-n)$ if $t < -n$. Let $u_n \in H_0^1(\Omega)$ be the
solution of (B.40) corresponding to $g_n$.
Note that $u_n$ satisfies
$$
\int_\Omega \nabla u_n \cdot \nabla v + \int_\Omega g_n(u_n) v =
\langle T, v \rangle \quad \forall v \in H_0^1(\Omega).
$$
Using $v = u_n$ as a test function, we get
$$
\int_\Omega |\nabla u_n|^2 + \int_\Omega g_n(u_n) u_n = \langle T,
u_n \rangle \leq C \left( \int_\Omega |\nabla u_n|^2 \right)^{1/2}.
$$
Thus,
$$
\int_\Omega g_n(u_n) u_n \leq C \quad \text{and} \quad \int_\Omega
|\nabla u_n|^2 \leq C, \tag B.41
$$
for some constant $C > 0$ independent of $n \geq 1$.
Since $(u_n)$ is uniformly bounded in $H_0^1(\Omega)$, then up to a
subsequence we can find $u \in H_0^1(\Omega)$ such that
$$
u_n \to u \quad \text{in $L^1$ and a.e.}
$$
By (B.41), for any $M > 0$, we also have
$$
\underset{[|u_n| \geq M]}\to\int |g_n(u_n)| \leq \frac{1}{M} \int_\Omega g_n(u_n)
u_n \leq \frac{C}{M}.
$$
We claim that
$$
g_n(u_n) \quad \text{is equi-integrable.}
$$
In fact, for any Borel set $E \subset \Omega$, we estimate
$$
\int_E |g_n(u_n)| = \underset{[|u_n| < M]}\to{\int_E} |g_n(u_n)| +
\underset{[|u_n| \geq M]}\to{\int_E} |g_n(u_n)| \leq A_M |E| + \frac{C}{M},
$$
where $A_M = \max{\{g(M), -g(-M)\}}$. Given $\eps > 0$, let $M > 0$
sufficiently large so that $\frac{C}{M}
< \eps$. With $M$ fixed, we take $|E|$ small enough so that $ A_M |E|
< \eps$. We conclude that
$$
\int_E |g_n(u_n)| < 2 \eps \quad \forall n \geq 1.
$$
Thus, $(g_n(u_n))$ is equi-integrable. Since $u_n \to u$ a.e., it follows
from Egorov's lemma that $g_n(u_n) \to g(u)$ in
$L^1(\Omega)$. Therefore, $u$ satisfies (B.40). By Proposition~B.3, this
solution is unique.
\enddemo

\smallskip
Combining the techniques from both proof, we have the following:

\proclaim{Theorem B.4}
For every $f \in L^1(\Omega)$ and every $T \in H^{-1}(\Omega)$, the equation
$$
\left\{ \alignedat2
-\Delta u + g(u)  & = f + T  && \quad \text{in }  \Omega,\\
u & = 0 && \quad \text{on } \partial \Omega,
\endalignedat \right.\tag B.42
$$
has a unique solution $u \in L^1(\Omega)$ with $g(u) \in L^1(\Omega)$.
\endproclaim

\demo{Proof}
Let $f_n$ be a sequence in $C^\infty(\overline\Omega)$ converging to
$f$ in $L^1(\Omega)$. Since $f_n + T \in H^{-1}$, we can apply Theorem~B.3
to obtain a solution $u_n$ of (B.42) for $f_n + T$. For every $n_1,n_2
\geq 1$, we have
$$
-\Delta (u_{n_1} - u_{n_2}) + g(u_{n_1}) - g(u_{n_2}) = f_{n_1} -
 f_{n_2} \quad \text{in } (C_0^2)^*. \tag B.43
$$
It follows from Proposition~B.3 that
$$
\int_\Omega \big| g(u_{n_1}) - g(u_{n_2}) \big| \leq \int_\Omega |f_{n_1} - f_{n_2}|.
$$
Thus, $(g(u_n))$ is a Cauchy sequence. Returning to (B.43), we
conclude from (B.3) that $(u_n)$ is
Cauchy in $L^1(\Omega)$. Passing to the limit as $n \to \infty$, we
find a solution $u \in L^1(\Omega)$ of (B.42). By Proposition~B.3, the
solution is unique.
\enddemo

\proclaim{Corollary B.3}
Let $\mu \in \Cal M(\Omega)$. If $\mu$ is diffuse, then \rom{(B.17)} admits
a unique solution $u \in L^1(\Omega)$ with $g(u) \in L^1(\Omega)$.
\endproclaim

\demo{Proof}
It suffices to observe that, by a result of
Boccardo-Gallou\"et-Orsina~[BGO1], every diffuse measure
$\mu$ belongs to $L^1 + H^{-1}$.
\enddemo

\smallskip
Concerning the existence of solutions for {\it every}\/ measure $\mu
\in \Cal M(\Omega)$, we have

\proclaim{Theorem B.5 (B\'enilan-Brezis [BB])}
Assume $N \geq 2$ and
$$
|g(t)| \leq C (|t|^p + 1) \quad \forall t \in \Bbb R,\tag B.44
$$
for some $p < \frac{N}{N-2}$. Then, for every $\mu \in \Cal
M(\Omega)$, problem \rom{(B.17)} has a unique solution $u \in
L^1(\Omega)$.
\endproclaim

\smallskip
Assumption (B.44) is optimal, in the sense that if $N \geq 3$, $g(t) = |t|^{p-1}t$
and $p \geq \frac{N}{N-2}$, then (B.17) has no weak solution for $\mu = \delta_a$,
where $a \in \Omega$:

\proclaim{Theorem B.6 (B\'enilan-Brezis [BB]; Brezis-V\'eron [BV])}
Assume $N \geq 3$. If $p \geq \frac{N}{N-2}$, then, for any $a \in
\Omega$, the problem
$$
\left\{
\alignedat2
-\Delta u + |u|^{p-1}u & = \delta_a && \quad \text{in }  \Omega,\\
u & = 0  && \quad \text{on } \partial \Omega,
\endalignedat \right.
$$
has no solution $u \in L^p(\Omega)$.
\endproclaim

\bigskip


\subhead Appendix C: Correspondence between $\big[ \bs C_{\pmb 0}
\pmb(\overline{\bs\Omega} \pmb)\big]^{\pmb *}$
and $\big[ \bs C \pmb(\overline{\bs\Omega} \pmb) \big]^{\pmb *}$  \endsubhead
\medskip

In this section we establish the following

\proclaim{Proposition C.1}
Given $\mu \in \big[ C_0(\overline\Omega)\big]^*$, there exists a
unique $\tilde \mu \in \big[ C(\overline\Omega)\big]^*$ such that
$$
\tilde \mu = \mu \quad \text{on } C_0(\overline\Omega) \quad
\text{and} \quad |\tilde \mu|(\partial\Omega) = 0. \tag C.1
$$
In addition, the map $\mu \mapsto \tilde \mu$ is a linear isometry.
\endproclaim

In order to prove Proposition C.1, we shall need the following

\proclaim{Lemma C.1}
Given $\varep > 0$, there exists $\delta > 0$ such that if $\zeta \in
C_0(\overline\Omega)$, $|\zeta| \leq 1$ in $\overline\Omega$, and
$\supp{\zeta} \subset \overline\Omega \backslash \Omega_\delta$, then
$$
|\langle \mu , \zeta \rangle| \leq \varep.
$$
\endproclaim

Here, we denote by $\Omega_\delta$ the set $\big\{ x \in \Omega \, ; \, d(x,
\partial\Omega)> \delta \big\}$.

\demo{Proof}
We argue by contradiction. Assume there exist $\varep_0 > 0$ and a
sequence $(\zeta_n) \subset C_0(\overline\Omega)$ such that $|\zeta_n|
\leq 1$ in $\overline\Omega$, $\supp{\zeta_n} \subset \overline\Omega
\backslash \Omega_{1/n}$, and
$$
\langle \mu, \zeta_n \rangle > \varep_0 \quad \forall n \geq 1.
$$
Without loss of generality, we may assume that each $\zeta_n$ has
compact support in $\Omega$ (this is always possible, by density of
$C_\trc^\infty(\Omega)$ in $C_0(\overline\Omega)$). In
particular, we can extract a subsequence $(\zeta_{n_j})$ such that
$\supp{\zeta_{n_j}}$ are all disjoint. For any $k \geq 1$, let $\tilde
\zeta_k = \sum_{j=1}^k{\zeta_{n_j}}$. By construction,
$$
\| \tilde \zeta_k \|_{L^\infty} \leq 1 \quad \text{and} \quad \supp{\tilde \zeta_k} \subset \Omega.
$$
Moreover,
$$
k \varep_0 < \langle \mu, \tilde\zeta_k \rangle \leq \|\mu\|_{\Cal M}.
$$
Since $k \geq 1$ was arbitrary, this gives a contradiction.
\enddemo

\medskip

\demo{Proof of Proposition C.1}
Let $\mu \in \big[C_0(\overline\Omega)\big]^*$. Given $\zeta \in
C(\overline\Omega)$, let $(\zeta_n)$ be any sequence in
$C_0(\overline\Omega)$ such that
$$
\| \zeta_n \|_{L^\infty} \leq C \quad \text{and} \quad \zeta_n \to
\zeta \quad \text{in $L^\infty_{\text{\rm loc}}(\Omega)$}.
$$
It easily follows from Lemma~C.1 that $(\langle \mu, \zeta_n \rangle)$
is Cauchy in $\Bbb R$. In particular, the limit $\displaystyle \lim_{n \to
\infty}{\langle \mu, \zeta_n \rangle}$ exists and is independent of
the sequence $(\zeta_n)$. Set
$$
\langle \tilde \mu, \zeta \rangle = \lim_{n \to \infty}{\langle \mu, \zeta_n \rangle}.
$$
Clearly, $\tilde\mu$ is a continuous linear functional on
$C(\overline\Omega)$ and
$$
\langle \tilde \mu, \zeta \rangle = \langle \mu, \zeta \rangle \quad
\forall \zeta \in C_0(\overline\Omega).
$$
In addition, Lemma~C.1 implies that $|\tilde
\mu|(\partial\Omega) = 0$; in particular, $\|\tilde\mu\|_{C^*} =
\|\mu\|_{(C_0)^*}$. The uniqueness of $\tilde\mu$ follows immediately from (C.1).
\enddemo

\bigskip


\subhead Appendix D: A new decomposition for diffuse measures  \endsubhead
\medskip

The goal of this section is to establish Theorem~3. Let $G$ denote the
Green function of the Laplacian in $\Omega$. Given $\mu \in \Cal M(\Omega)$, $\mu \geq
0$, set
$$
G(\mu)(x) = \int_\Omega G(x,y) \, d\mu(y).
$$
Note that $G(\mu)$ is well-defined for every $x \in \Omega$, possibly
taking values $+\infty$.

We first present some well-known results in Potential Theory:

\proclaim{Lemma D.1}
Let $\mu \in \Cal M(\Omega)$, $\mu \geq 0$, be such that $G(\mu) <
\infty$ everywhere in $\Omega$. Given $\eps > 0$, there exists $L
\subset \Omega$ compact such that
$$
\mu(\Omega \backslash L) < \eps \quad \text{and} \quad G(\mu\lfloor_L)
\in C_0(\overline\Omega). \tag D.1
$$
\endproclaim

\demo{Proof}
If $\mu$ has compact support in $\Omega$, then Lemma~D.1 is precisely
Theorem~6.21 in [H]. For an arbitrary $\mu \in \Cal M(\Omega)$, $\mu
\geq 0$, such that $G(\mu) < \infty$ in $\Omega$, we proceed as
follows. By inner regularity of $\mu$, there
exists $K \subset \Omega$ compact such that $\mu(\Omega\backslash K) <
\frac{\eps}{2}$. Since $G(\mu\lfloor_K) \leq G(\mu)$, the function
$G(\mu\lfloor_K)$ is also finite everywhere in $\Omega$. Then, by Theorem~6.21 in
[H], there exists $L \subset \Omega$ compact such that
$$
\mu\lfloor_K(\Omega \backslash L) < \frac{\eps}{2} \quad \text{and} \quad
G(\mu\lfloor_{K \cap L}) \in C_0(\overline\Omega).
$$
We conclude that (D.1) holds with $L$ replaced by $K \cap L$.
\enddemo

As a consequence of Lemma~D.1, we have

\proclaim{Proposition D.1}
Let $u \in W^{1,1}_0(\Omega)$ be such that $\Delta u$ is a diffuse measure in
$\Omega$. Then, there exists a sequence $(u_n) \subset
C_0(\overline\Omega)$ such that $\Delta u_n \in \Cal
M(\Omega)$, $\forall n \geq 1$,
$$
u = \sum_{n=1}^\infty{u_n} \quad \text{a.e.$ $ in $\Omega$} \quad \text{and} \quad
\|\Delta u\|_{\Cal M} = \sum_{n=1}^\infty{\|\Delta u_n\|_{\Cal
M}}. \tag D.2
$$
\endproclaim

\demo{Proof}
We shall split the proof of Proposition~D.1 into three steps.
\smallskip

\nd
{\it Step 1.} Let $\mu \geq 0$ be a measure such that
$G(\mu) < \infty$ everywhere in $\Omega$. Then, there exist disjoint
Borel sets $A_n \subset \Omega$ such that
$$
\mu\Big( \Omega \backslash \bigcup_{n=1}^\infty{A_n} \Big) = 0 \quad
\text{and} \quad G(\mu\lfloor_{A_n}) \in C_0(\overline\Omega) \quad
\forall n \geq 1. \tag D.3
$$

This result easily follows from Lemma~D.1 by an induction argument.

\medskip
\nd
{\it Step 2.} Let $\mu \geq 0$ be a diffuse measure in $\Omega$. Then, there exist
disjoint Borel sets $A_n \subset \Omega$ such that
$$
\mu\Big( \Omega \backslash \bigcup_{n=1}^\infty{A_n} \Big) = 0 \quad
\text{and} \quad G(\mu\lfloor_{A_n}) \in C_0(\overline\Omega) \quad
\forall n \geq 1. \tag D.4
$$

For each $k \geq 1$, let
$$
E_k = \big\{ x\in \Omega \; ;\; G(\mu)(x) \leq k  \big\}.
$$
Since $G(\mu)$ is lower semicontinuous (by Fatou), $E_k$ is
closed in $\Omega$. Clearly, we have $G(\mu\lfloor_{E_k}) \leq k$ in $E_k$
and $G(\mu\lfloor_{E_k})$ is harmonic in $\Omega \backslash E_k$.
Therefore, by the maximum principle, $G(\mu\lfloor_{E_k}) \leq k$
everywhere in $\Omega$.

\nd
Applying the previous step to the measures
$\mu\lfloor_{E_k \backslash E_{k-1}}$, one can find disjoint Borel
sets $A_n \subset \Omega$ such that
$$
\mu\Big( F \backslash \bigcup_{n=1}^\infty{A_n} \Big) = 0 \quad
\text{and} \quad G(\mu\lfloor_{A_n}) \in C_0(\overline\Omega) \quad
\forall n \geq 1,
$$
where
$$
F = \big\{ x\in \Omega \; ;\; G(\mu)(x) < \infty \big\}.
$$
Since $\mu$ is diffuse and $\Omega \backslash F$ has zero capacity
(see e.g. [H, Theorem~7.33]), we have $\mu(\Omega \backslash F) =
0$. Thus,
$$
\mu\Big( \Omega \backslash \bigcup_{n=1}^\infty{A_n} \Big) = 0,
$$
from which the result follows.

\medskip
\nd
{\it Step 3.} Proof of Proposition~D.1 completed.
\smallskip

Set $\mu = - \Delta u$. Applying Step~2 to $\mu^+$, one can find
disjoint Borel sets $(A_n)$ such that
$$
\mu^+\Big( \Omega \backslash \bigcup_{n=1}^\infty{A_n} \Big) = 0 \quad
\text{and} \quad G(\mu^+\lfloor_{A_n}) \in C_0(\overline\Omega) \quad
\forall n \geq 1.
$$
Similarly, there exist disjoint Borel sets $(B_n)$ such that
$$
\mu^-\Big( \Omega \backslash \bigcup_{n=1}^\infty{B_n} \Big) = 0 \quad
\text{and} \quad G(\mu^-\lfloor_{B_n}) \in C_0(\overline\Omega) \quad
\forall n \geq 1.
$$
Since
$$
\mu = \mu^+ - \mu^- = \sum_{n=1}^\infty{\mu^+\lfloor_{A_n}} -
\sum_{n=1}^\infty{\mu^-\lfloor_{B_n}},
$$
we have
$$
u = \sum_{n=1}^\infty{G(\mu^+\lfloor_{A_n})} -
\sum_{n=1}^\infty{G(\mu^-\lfloor_{B_n})} \quad \text{a.e.}
$$
and
$$
\|\Delta u\|_{\Cal M} = \sum_{n=1}^\infty{\big\|\mu^+\lfloor_{A_n}\big\|_{\Cal
M}} + \sum_{n=1}^\infty{\big\|\mu^-\lfloor_{B_n}\big\|_{\Cal M}}.
$$
This concludes the proof of the proposition.
\enddemo

\medskip
We can now present the

\demo{Proof of Theorem 3}
Let $u \in W^{1,1}_0(\Omega)$ be the unique solution of
$$
- \Delta u = \mu \quad \text{in $(C_0^2)^*$.}
$$
Let $(u_n) \subset C_0(\overline\Omega)$ be the
sequence given by Proposition~D.1.
For $\delta>0$ fixed, take $w_n \in
C_0^2(\overline\Omega)$ such that
$$
\|u_n - w_n\|_{L^\infty} \leq \frac{\delta}{2^n} \quad
\text{and}
\quad
\|\Delta w_n\|_{L^1} \leq \|\Delta u_n\|_{\Cal M}.
$$
Let
$$
v = \sum_{n=1}^\infty{ (u_n - w_n) } \quad
\text{and} \quad f = - \sum_{n=1}^{\infty}{\Delta w_n}.
$$
Since
$$
\| v \|_{L^\infty} \leq \sum_{n=1}^\infty{\|u_n - w_n \|_{L^\infty}}
\leq \delta, \tag D.5
$$
we have $v \in C_0(\overline\Omega)$ and $\|v\|_{L^\infty} \leq \delta$.
Moreover,
$$
\|f\|_{L^1} \leq \sum_{n=1}^{\infty}{\|\Delta w_n \|_{L^1}} \leq
\sum_{n=1}^{\infty}\| \Delta u_n\|_{\Cal M} = \| \mu \|_{\Cal M} \tag D.6
$$
implies $f \in L^1(\Omega)$. Finally, by construction, we have
$$
\mu = f - \Delta v \quad \text{in } (C_0^2)^*. \tag D.7
$$
In particular, $\Delta v = f - \mu$ is a measure and $\|\Delta
v\|_{\Cal M} \leq 2 \| \mu \|_{\Cal M}$. Thus,
$$
\|\nabla v\|_{L^2}^2 \leq \|v\|_{L^\infty}\| \Delta v\|_{\Cal M} \leq
2 \delta \|\mu\|_{\Cal M}. \tag D.8
$$
Since $v \in C_0(\overline\Omega) \cap H_0^1$, (0.21) immediately follows from
(D.7). Moreover, replacing $\delta$ by $\frac{\delta}{2} \|\mu\|_{\Cal
M}$ in (D.5) and (D.8), we conclude that (0.22) holds. The proof of
Theorem~3 is complete.
\enddemo

\medskip

Note that our construction of $f \in L^1$ and $v \in L^\infty$ satisfying
(0.21) is not linear with respect to $\mu$. Here is a natural question:

\definition{Open problem 7}
Can one find a bounded linear operator
$$
T \; : \; \mu \in \Cal M_\trd(\Omega) \; \longmapsto \; (f,v) \in L^1 \times L^\infty
$$
such that (0.21) and (0.22) hold?
\enddefinition

After receiving a preprint of our work, A. Ancona~[A2] has provided a
negative answer to the question above.

\bigskip


\subhead Appendix E: Equivalence between $\text{cap}_{\bs H^{\pmb 1}}$
and $\text{cap}_{\bs \Delta, \pmb 1}$
\endsubhead
\medskip

Given a compact set $K \subset \Omega$, let $\capt_{\Delta,1}{(K)}$ denote
the capacity associated to the Laplacian. More precisely,
$$
\capt_{\Delta,1}(K) = \inf{\left\{ \int_\Omega |\Delta \varphi|\; ;
\; \varphi \in C_\trc^\infty(\Omega), \; \varphi \geq 1 \text{
in some neighborhood of } K \right\}}.
$$
In order to avoid confusion, throughout this section we shall denote
by $\capt_{H^1}$ the Newtonian capacity with respect to $\Omega$
(which we simply denote $\capt$ everywhere else in this paper).

\smallskip
The main result in this appendix is the following

\proclaim{Theorem E.1}
For every compact set $K \subset \Omega$, we have
$$
\capt_{\Delta,1}(K) = 2 \capt_{H^1} (K). \tag E.1
$$
\endproclaim

\remark{Remark E.1}
In an earlier version of this work, we had only established the
equivalence between $\capt_{H^1}$ and $\capt_{\Delta,1}$. The exact
formula (E.1) has been suggested to us by A. Ancona.
\endremark

\medskip

We first prove the following

\proclaim{Lemma E.1}
Let $K \subset \Omega$ be a compact set. Given $\eps > 0$, there
exists $\psi \in C_\trc^\infty(\Omega)$ such that $0 \leq \psi \leq 1$
in $\Omega$, $\psi = 1$ in some neighborhood of $K$, and
$$
\int_\Omega |\Delta \psi| \leq 2 \capt_{H^1} (K) + \eps. \tag E.2
$$
\endproclaim

\demo{Proof}
Let $\omega \subset\subset \Omega$ be an open set such that $K \subset
\omega$ and
$$
\capt_{H^1}(\overline\omega) \leq \capt_{H^1}(K) + \frac{\eps}{4}.
$$
Let $u$ denote the capacitary potential of $\overline\omega$. More
precisely, let $u \in H^1_0(\Omega)$ be such that $u = 1$ in
$\overline\omega$ and
$$
\int_\Omega |\nabla u|^2 = \capt_{H^1}(\overline\omega).
$$
Note that $u$ is superharmonic in $\Omega$ and harmonic in $\Omega
\backslash \overline\omega$. In particular, $0 \leq u \leq 1$. Since
$\supp{\Delta u} \subset [u=1]$, $u$ is continuous (see [H,
Theorem~6.20]) and
$$
\|\Delta u\|_{\Cal M} = - \int_\Omega \Delta u = - \int_\Omega u \Delta
u = \int_\Omega |\nabla
u|^2 = \capt_{H^1}(\overline\omega).
$$
Given $\delta > 0$ small, set
$$
v = \frac{(u-\delta)^+}{1-\delta}.
$$
Since $v$ has compact support in $\Omega$, we have
$$
\int_\Omega \Delta v = 0. \tag E.3
$$
Moreover, $\Delta v$ is a diffuse measure (note that $v \in
H^1_0(\Omega)$) and
$$
\supp{\Delta v} \subset [v = 0] \cup [v = 1]. \tag E.4
$$
Thus, by Corollary~1.3 in [BP2], we have
$$
\Delta v \geq 0 \quad \text{in } [v = 0] \quad \text{and} \quad \Delta v
\leq 0 \quad \text{in } [v = 1]. \tag E.5
$$
It then follows from (E.3)--(E.5) that
$$
\| \Delta v \|_{\Cal M}  = 2 \underset{[v = 1]}\to\int |\Delta v|.
$$
Since $\Delta v = \frac{1}{1 - \delta} \Delta u$ in $[v = 1]$, we
conclude that
$$
\| \Delta v \|_{\Cal M} \leq \frac{2}{1 - \delta} \|\Delta u\|_{\Cal M}.
$$
Using the same notation as in Section~4, we now take $n \geq 1$
sufficiently large so that the function $\psi = \rho_n * v$ has
compact support in $\Omega$ and $\psi = 1$ in some neighborhood of
$K$. We claim that $\psi$ satisfies all the required properties.
In fact, since $0 \leq \psi \leq 1$ in $\Omega$, we only have to show
that (E.2) holds. Note that
$$
\int_\Omega |\Delta \psi| \leq \|\Delta v\|_{\Cal M} \leq \frac{2}{1 -
\delta} \|\Delta u\|_{\Cal M} = \frac{2}{1 - \delta} \capt_{H^1}(\overline\omega).
$$
Choosing $\delta > 0$ so that
$$
\frac{\delta}{ 1 - \delta} \capt_{H^1}(\overline\omega) <
\frac{\eps}{4},
$$
we have
$$
\int_\Omega |\Delta \psi| \leq 2 \left( 1 + \frac{\delta}{ 1 - \delta}
\right)\capt_{H^1}(\overline\omega) \leq 2 \capt_{H^1}(K) + \eps,
$$
which is precisely (E.2).
\enddemo

We now present the

\demo{Proof of Theorem~E.1}
In view of Lemma~E.1, it suffices to show that
$$
\capt_{H^1}(K) \leq \frac{1}{2} \capt_{\Delta,1}(K). \tag E.6
$$
Let $\varphi \in C_\trc^\infty(\Omega)$ be such that $\varphi \geq 1$ in
some neighborhood of $K$. Set
$\tilde\varphi = \min{\{ 1, \varphi^+\}}$. For $n \geq 1$ sufficiently
large, the function $\tilde \varphi_n = \rho_n * \tilde\varphi$ belongs to
$C_\trc^\infty(\Omega)$ and $\tilde \varphi_n = 1$ in some neighborhood of $K$. We
then have
$$
\capt_{H^1}(K) \leq \int_\Omega |\nabla \tilde\varphi_n|^2 \leq
\int_\Omega |\nabla \tilde\varphi|^2 = \int_\Omega \nabla
\tilde\varphi \cdot \nabla \varphi = - \int_\Omega \tilde\varphi
\Delta \varphi.
$$
Recall that $\varphi$ has compact support in $\Omega$ and $0 \leq
\tilde \varphi \leq 1$. Thus, $\int_\Omega
\Delta\varphi = 0$ and we have
$$
\capt_{H^1}(K) \leq - \int_\Omega \Big( \tilde\varphi  - \frac{1}{2}
\Big) \Delta \varphi \leq \frac{1}{2} \int_\Omega |\Delta\varphi|.
$$
Since $\varphi$ was arbitrary, we conclude that (E.6) holds. This
establishes Theorem~E.1.
\enddemo

\bigskip


\noindent{\it Acknowledgments:\/} We warmly thank A. Ancona for
enlightening discussions and suggestions. The first author (H.B.) and
the second author (M.M.) are partially sponsored by an E.C. Grant
through the RTN Program ``Front-Singularities'', HPRN-CT-2002-00274.
H.B. is also a member of the Institut Universitaire de France.  



\bigskip
\bigskip
\Refs
\widestnumber\key{AAAAA}


\ref\key A1
\by A. Ancona
\paper Une propri\'et\'e d'invariance des ensembles absorbants par
              perturbation d'un op\'erateur elliptique
\jour Comm. Partial Differential Equations
\vol 4
\yr 1979
\pages 321--337
\endref

\ref\key A2
\bysame
\paper Sur une question de H.~Brezis, M.~Marcus et A.C.~Ponce
\jour Ann. Inst. H. Poincar\'e Anal. Non Lin\'eaire
\vol 23
\yr 2006
\pages 127--133
\endref

\ref\key BP
\by P. Baras and M. Pierre
\paper Singularit\'es \'eliminables pour des \'equations
              semi-lin\'eaires
\jour Ann. Inst. Fourier (Grenoble)
\vol 34
\yr 1984
\pages 185--206
\endref

\ref\key BLOP
\by D. Bartolucci, F. Leoni, L. Orsina and A.C. Ponce
\paper Semilinear equations with exponential nonlinearity and measure data
\jour Ann. Inst. H. Poincar\'e Anal. Non Lin\'eaire
\vol 22
\yr 2005
\pages 799--815
\endref

\ref\key Ba
\by J.R. Baxter
\paper Inequalities for potentials of particle systems
\jour Illinois J. Math.
\vol 24
\yr 1980
\pages 645--652
\endref

\ref\key BM
\by G. Barles and F. Murat
\paper Uniqueness and the maximum principle for quasilinear elliptic
              equations with quadratic growth conditions
\jour Arch. Rational Mech. Anal.
\vol 133
\yr 1995
\pages 77--101
\endref

\ref\key BB
\by Ph. B\'enilan and H. Brezis
\paper Nonlinear problems related to the Thomas-Fermi equation
\jour J. Evol. Equ.
\vol 3
\yr 2004
\pages 673--770. Dedicated to Ph.~B\'enilan
\endref

\ref\key BGO1
\by L. Boccardo, T. Gallou\"et and L. Orsina
\paper Existence and uniqueness of entropy solutions for nonlinear
              elliptic equations with measure data
\jour Ann. Inst. H. Poincar\'e Anal. Non Lin\'eaire
\vol 13
\yr 1996
\pages 539--551
\endref

\ref\key BGO2
\bysame
\paper Existence and nonexistence of solutions for some nonlinear
              elliptic equations
\jour J. Anal. Math.
\vol 73
\yr 1997
\pages 203--223
\endref

\ref\key B1
\by H. Brezis
\book Nonlinear problems related to the Thomas-Fermi
equation. \rom{In: Contemporary developments in
continuum mechanics and partial differential equations (G.M. de la
Penha and L.A. Medeiros, eds.) Proc. Internat. Sympos., Inst. Mat.,
Univ. Fed. Rio de Janeiro, Rio de Janeiro}
\publ North Holland
\publaddr Amsterdam
\yr 1978, pp. 74--80
\endref

\ref\key B2
\bysame
\book Some variational problems of the {T}homas-{F}ermi type. \rom{In:
Variational inequalities and complementarity problems (R.W. Cottle,
F. Giannessi and J.-L. Lions, eds.) Proc. Internat. School, Erice,
1978}
\publ Wiley
\publaddr Chichester
\yr 1980, pp. 53--73
\endref

\ref\key B3
\bysame
\book Probl\`emes elliptiques et paraboliques non lin\'eaires avec
donn\'ees mesures. \rom{Gou\-laouic-Meyer-Schwartz Seminar, 1981/1982}
\publ \'Ecole Polytech.
\publaddr Palaiseau
\yr 1982, pp. X.1--X.12
\endref

\ref\key B4
\bysame
\book Nonlinear elliptic equations involving measures. \rom{In:
Contributions to nonlinear partial differential equations (C.~Bardos,
A.~Damlamian, J.I.~Diaz and J.~Hernandez, eds.) Madrid, 1981}
\publ Pitman
\publaddr Boston, MA
\yr 1983, pp. 82--89
\endref

\ref\key B5
\bysame
\paper Semilinear equations in {${\Bbb R}\sp N$} without condition at infinity
\jour Appl. Math. Optim.
\vol 12
\yr 1984
\pages 271--282
\endref

\ref\key BBr
\by H. Brezis and F.E. Browder
\paper Strongly nonlinear elliptic boundary value problems
\jour Ann. Scuola Norm. Sup. Pisa Cl. Sci.
\vol 5
\yr 1978
\pages 587--603
\endref

\ref\key BCMR
\by H. Brezis, T. Cazenave, Y. Martel and A. Ramiandrisoa
\paper Blow up for {$u\sb t-\Delta u=g(u)$} re\-vis\-ited
\jour Adv. Differential Equations
\vol 1
\yr 1996
\pages 73--90
\endref

\ref\key BMP
\by H. Brezis, M. Marcus and A.C. Ponce
\paper A new concept of reduced measure for nonlinear elliptic equations
\jour C. R. Acad. Sci. Paris, Ser. I
\vol 339
\yr 2004
\pages 169--174
\endref

\ref\key BP1
\by H. Brezis and A.C. Ponce
\paper Remarks on the strong maximum principle
\jour Differential Integral Equations
\vol 16
\yr 2003
\pages 1--12
\endref

\ref\key BP2
\bysame
\paper Kato's inequality when $\Delta u$ is a measure
\jour C. R. Acad. Sci. Paris, Ser. I
\vol 338
\yr 2004
\pages 599--604
\endref

\ref\key BP3
\bysame
\paper Reduced measures on the boundary
\jour J. Funct. Anal
\vol 229
\yr 2005
\pages 95--120
\endref

\ref\key BP4
\bysame
\paper Reduced measures for obstacle problems
\jour Adv. Diff. Equations
\vol 10
\yr 2005
\pages 1201--1234
\endref

\ref\key BSe
\by H. Brezis and S. Serfaty
\paper A variational formulation for the two-sided obstacle problem
with measure data
\jour Commun. Contemp. Math.
\vol 4
\yr 2002
\pages 357--374
\endref

\ref\key BS
\by H. Brezis and W.A. Strauss
\paper Semilinear second-order elliptic equations in {$L\sp{1}$}
\jour J. Math. Soc. Japan
\vol 25
\yr 1973
\pages 565--590
\endref

\ref\key BV
\by H. Brezis and L. V\'eron
\paper Removable singularities for some nonlinear elliptic equations
\jour Arch. Rational Mech. Anal.
\vol 75
\yr 1980/81
\pages 1--6
\endref

\ref\key DD
\by P. Dall'Aglio and G. Dal Maso
\paper Some properties of the solutions of obstacle problems with
              measure data
\jour Ricerche Mat.
\vol 48
\yr 1999
\pages suppl., 99--116. Papers in memory of Ennio De Giorgi
\endref

\ref\key DM
\by C. Dellacherie and P.-A. Meyer
\book Probabilit{\'e}s et potentiel, \rom{Chapitres I {\`a} IV,
              Publications de l'Institut de Math{\'e}matique de l'Universit{\'e}
              de Strasbourg, No. XV,
              Actualit{\'e}s Scientifiques et Industrielles, No. 1372}
\publ Hermann
\publaddr Paris
\yr 1975
\endref

\ref\key DS
\by N. Dunford and J.T. Schwartz
\book Linear operators. {P}art {I}
\publ Wiley
\publaddr New York
\yr 1958
\endref

\ref\key DP
\by L. Dupaigne and A.C. Ponce
\paper Singularities of positive supersolutions in elliptic PDEs
\jour Selecta Math. (N.S.)
\vol 10
\yr 2004
\pages 341--358
\endref

\ref\key DVP
\by C. De La Vall\'ee-Poussin
\paper Sur l'int\'egrale de Lebesgue
\jour Trans. Amer. Math. Soc.
\vol 16
\yr 1915
\pages 435--501
\endref

\ref\key DPP
\by L. Dupaigne, A.C. Ponce and A. Porretta
\paper Elliptic equations with vertical asymptotes in the nonlinear term
\jour to appear
\endref

\ref\key D1
\by E.B. Dynkin
\book Diffusions, superdiffusions and partial differential
              equations
\publ American Mathematical Society
\publaddr Providence, RI
\yr 2002
\endref

\ref\key D2
\bysame
\book Superdiffusions and positive solutions of nonlinear partial differential equations
\publ American Mathematical Society
\publaddr Providence, RI
\yr 2004
\endref

\ref\key FTS
\by M. Fukushima, K. Sato and S. Taniguchi
\paper On the closable parts of pre-{D}irichlet forms and the fine
              supports of underlying measures
\jour Osaka J. Math.
\vol 28
\yr 1991
\pages 517--535
\endref

\ref\key GM
\by T. Gallou{\"e}t and J.-M. Morel
\paper Resolution of a semilinear equation in {$L\sp{1}$}
\jour Proc. Roy. Soc. Edinburgh Sect. A
\vol 96
\yr 1984
\pages 275--288
\moreref
\jour Corrigenda: Proc. Roy. Soc. Edinburgh Sect. A {\bf 99} (1985), 399
\endref

\ref\key GV
\by A. Gmira and L. V{\'e}ron
\paper Boundary singularities of solutions of some nonlinear elliptic
              equations
\jour Duke Math. J.
\vol 64
\yr 1991
\pages 271--324
\endref

\ref\key GV
\by M. Grillot and L. V{\'e}ron
\paper Boundary trace of the solutions of the prescribed Gaussian
curvature equation.
\jour  Proc. Roy. Soc. Edinburgh Sect. A
\vol 130
\yr 2000
\pages 527--560
\endref


\ref\key GRe
\by M. Grun-Rehomme
\paper Caract\'erisation du sous-diff\'erential d'int\'egrandes
              convexes dans les espaces de {S}obolev
\jour J. Math. Pures Appl.
\vol 56
\yr 1977
\pages 149--156
\endref

\ref\key H
\by L.L. Helms
\book Introduction to potential theory
\publ Wiley-Interscience
\publaddr New York
\yr 1969
\endref

\ref\key K
\by T. Kato
\paper Schr\"odinger operators with singular potentials
\jour Israel J. Math.
\vol 13
\yr 1972
\pages 135--148 (1973)
\endref

\ref\key LG1
\by J.-F. Le Gall
\paper The {B}rownian snake and solutions of {$\Delta u=u\sp 2$} in a
              domain
\jour Probab. Theory Related Fields
\vol 102
\yr 1995
\pages 393--432
\endref

\ref\key LG2
\bysame
\paper A probabilistic {P}oisson representation for positive
              solutions of {$\Delta u=u\sp 2$} in a planar domain
\jour Comm. Pure Appl. Math.
\vol 50
\yr 1997
\pages 69--103
\endref

\ref\key LS
\by E.H. Lieb and B. Simon
\paper  The Thomas-Fermi theory of atoms, molecules and solids
\jour  Advances in Math.
\vol 23
\yr 1977
\pages 22--116
\endref


\ref\key MV1
\by M. Marcus and L. V\'eron
\paper The boundary trace of positive solutions of semilinear
              elliptic equations: the subcritical case
\jour Arch. Rational Mech. Anal.
\vol 144
\yr 1998
\pages 201--231
\endref

\ref\key MV2
\bysame
\paper The boundary trace of positive solutions of semilinear
              elliptic equations: the supercritical case
\jour J. Math. Pures Appl.
\vol 77
\yr 1998
\pages 481--524
\endref

\ref\key MV3
\bysame
\paper Removable singularities and boundary traces
\jour J. Math. Pures Appl.
\vol 80
\yr 2001
\pages 879--900
\endref

\ref\key MV4
\bysame
\paper Capacitary estimates of solutions of a class of nonlinear
              elliptic equations
\jour C. R. Acad. Sci. Paris, Ser. I
\vol 336
\yr 2003
\pages 913--918
\endref

\ref\key MV5
\bysame
\paper Capacitary estimates of positive solutions of semilinear elliptic equations with absorption
\jour J. Eur. Math. Soc.
\vol 6
\yr 2004
\pages 483-527
\endref

\ref\key MV6
\bysame
\paper Nonlinear capacities associated to semilinear elliptic
equations
\jour in preparation
\endref


\ref\key P
\by A.C. Ponce
\book How to construct good measures. \rom{In: Elliptic and
parabolic problems (C. Bandle, H. Berestycki, B. Brighi, A. Brillard,
              M. Chipot, J.-M. Coron, C. Sbordone, I. Shafrir,
              V. Valente and G. Vergara-Caffarelli, eds.) Gaeta, 2004. A special tribute to the work of
Ha\"\i m Brezis}
\publ Birkh{\"a}user
\publaddr Basel, Boston, Berlin
\yr 2005, pp. 375--388
\endref

\ref\key Po
\by A. Porretta
\paper Absorption effects for some elliptic equations with
singularities
\jour Boll. Unione Mat. Ital. Sez. B Artic. Ric. Mat. (8)
\vol 8
\yr 2005
\pages 369--395
\endref

\ref\key S
\by G. Stampacchia
\book \'Equations elliptiques du second ordre \`a coefficients
discontinus
\publ Les Presses de l'Universit\'e de Montr\'eal
\publaddr Montr\'eal
\yr 1966
\endref

\ref\key Va
\by J.L. V\'azquez
\paper On a semilinear equation in {$\Bbb R\sp{2}$} involving bounded measures
\jour Proc. Roy. Soc. Edinburgh Sect. A
\vol 95
\yr 1983
\pages 181--202
\endref

\endRefs
\enddocument